\documentclass[10pt]{article}

\usepackage[utf8]{inputenc}
\usepackage[T1]{fontenc}
\usepackage{lmodern}
\usepackage[a4paper,left=2cm,right=2cm,top=2.75cm,bottom=2.75cm]{geometry}
\usepackage{amsmath,amssymb,amsfonts,amsthm,mathtools}
\usepackage{graphicx}
\usepackage{tikz}
\usepackage{placeins}
\usepackage{hyperref}
\usepackage{xcolor}
\hypersetup{hidelinks,bookmarksdepth=2}

\numberwithin{equation}{section}
\numberwithin{figure}{section}
\numberwithin{table}{section}

\newcommand{\eps}{\varepsilon}
\newcommand{\dd}{\,\mathrm d}
\newcommand{\TV}{\operatorname{TV}}
\newcommand{\loc}{\operatorname{loc}}
\newcommand{\diag}{\operatorname{diag}}
\newcommand{\cO}{\mathcal O}
\newcommand{\be}{\begin{equation}}
\newcommand{\ee}{\end{equation}}

\newcommand \bse {\begin{subequations}}
\newcommand \ese {\end{subequations}}

\newcommand \del \partial

\newcommand \Acal {\mathcal A}

\theoremstyle{plain}
\newtheorem{theorem}{Theorem}[section]
\newtheorem{proposition}[theorem]{Proposition}
\newtheorem{lemma}[theorem]{Lemma}
\newtheorem{corollary}[theorem]{Corollary}

\newcommand{\ang}[2]{\left\langle #1,#2\right\rangle}

\begin{document}

\title{The Riemann problem for a Lorentz-invariant model
\\
 of viscous relativistic fluids}

\author{Philippe G. LeFloch$^{1}$ and Iker Pascual-de-Zulueta$^{1,2}$}

\date{November 2026}

\maketitle

\footnotetext[1]{Laboratoire Jacques-Louis Lions,
Centre National de la Recherche Scientifique, Sorbonne Universit\'e,
4 Place Jussieu, 75252 Paris, France.
Email: {\tt contact@philippelefloch.org}.}

\footnotetext[2]{CEREMADE,
Universit\'e Paris Dauphine--PSL,
Place du Mar\'echal de Lattre de Tassigny, 75775 Paris Cedex 16, France.
\\
Email: {\tt iker.pascual-de-zulueta-barandiaran@dauphine.eu}.}

\begin{abstract}
We study the vanishing-viscosity limit for the one-dimensional relativistic Euler equations with large Riemann data and a Lorentz-invariant hyperbolic viscosity. For genuinely nonlinear pressure laws, including isothermal laws, with sound speed bounded away from zero and from the speed of light, we construct global self-similar profiles for arbitrary positive densities and finite rapidities, without a smallness restriction on the initial discontinuity. We establish positive lower and upper density bounds, rapidity bounds, and total-variation estimates, which are uniform as the viscosity parameter tends to zero. At positive viscosity the profiles equal the prescribed states outside the light cone and approach them super-exponentially in rapidity from inside the cone. For fixed data, the full family converges to the classical Euler Riemann solution, with differences tending to zero in $L^1$ on the whole rapidity line. The limit has bounded variation, takes the prescribed Riemann data, and satisfies every smooth convex entropy inequality. The proof combines weighted flux identities with a continuation of two half-profiles joined at a regular intermediate state. Near an isolated shock of either family, the centered profiles converge to a traveling wave that determines the leading width of the layer in the ray variable. Numerical experiments illustrate the layer asymptotics for large jumps and compare the model with spatial viscosity under Lorentz boosts.
\end{abstract}

%
%

\section{Introduction}
\label{section-1}

\subsection{Relativistic viscous fluids}

\paragraph{The PDE model of interest}
 
On a Lorentzian spacetime $(\mathcal M,g)$, a fluid of mass-energy density\footnote{Vacuum is avoided here by assuming the sound speed to be bounded away from zero.}  $\mu>0$, pressure $p(\mu)$, and future-directed unit velocity $U=(U^\alpha)$ has energy-momentum tensor\footnote{Here, $\alpha, \beta = 0,1,2,3,$ and Einstein's index convention is used.}
\bse
\be
\label{eq:energy-momentum-tensor}
T^{\alpha\beta}\coloneqq\bigl(\mu+p(\mu)\bigr)U^\alpha U^\beta+p(\mu)g^{\alpha\beta},\qquad g_{\alpha\beta}U^\alpha U^\beta=-1.
\ee
The pressure is prescribed as a function of the density. 
With $\nabla$ the Levi--Civita connection of the Lorentzian metric $g$, the \emph{relativistic Euler equations} express conservation of energy and momentum as
\be
\label{eq:relativistic-euler}
\nabla_\alpha T^{\alpha\beta}=0,
\ee
\ese
which is a first-order nonlinear hyperbolic system provided $p'(\mu)>0$. Solutions can involve shock waves and must be understood in the sense of distributions.  

We write $h(\mu)\coloneqq\mu+p(\mu)$ for the enthalpy density and define the sound speed by $c^2(\mu)\coloneqq p'(\mu)$, where the prime denotes differentiation with respect to $\mu$. Chaddha and Reintjes~\cite{ReintjesChaddha2026} recently introduced a viscosity term given by the wave operator in Minkowski spacetime. The corresponding \emph{relativistic Navier-Stokes equations} are the following covariant equations: 
\be
\label{eq:covariant-velocity-viscosity}
\nabla_\alpha \bigl( T^{\alpha\beta}(\mu,U) \bigr) =\nabla_\alpha\bigl(\eps\nabla^\alpha U^\beta\bigr).
\ee
Here the viscosity coefficient $\eps=\eps(\cdot,\mu,U)>0$ is a smooth scalar field defined on the spacetime that may depend on the spacetime point and on the fluid state $(\mu,U)$. The second-order part acts only on the velocity and is therefore \emph{degenerate in nature}, with derivatives of the density occurring only at first order. 

\paragraph{The viscous Riemann problem.}

Specifically, the setting above is simplified by working in two-dimensional Minkowski spacetime with metric $g=\diag(-1,1)$ and unit light speed. We also write $U=(\cosh\theta,\sinh\theta)$, where $\theta$ is the \emph{fluid rapidity}.
The Euler equations~\eqref{eq:relativistic-euler} read $\del_\alpha T^{\alpha\beta}=0$ in standard coordinates. We prescribe Riemann initial data for the density and rapidity, denoted by $z=(\mu,\theta)$:
\be
\label{eq:riemann-initial-data}
z(0,x)=\begin{cases}z_l,&x<0,\\z_r,&x>0,\end{cases}\qquad z_l=(\mu_l,\theta_l),\quad z_r=(\mu_r,\theta_r).
\ee
Moreover, we choose the viscosity coefficient $\eps=\eps(\cdot,\mu,U)>0$ so that the viscous approximation remains Lorentz-invariant. 

Hence, our main purpose here is 
 to establish \emph{uniform estimates} and study the singular-limit problem $\eps=\eps(\cdot,\mu,U) \to 0$  for the viscous Riemann problem in the setting of the one-dimensional relativistic Euler system with \emph{arbitrarily large data}; the endpoint density values $\mu_l$ and $\mu_r$ are arbitrary, and the endpoint velocities are any \emph{strictly subluminal velocities}, without any smallness requirement. For genuinely nonlinear pressure laws with uniformly positive, subluminal sound speed, our estimates indeed control the full family of self-similar solutions as the viscosity parameter tends to zero and remain uniform when the prescribed states range over a compact subset of the state space.

For a class of conservation laws, Dafermos~\cite{Dafermos_1973} and Tupciev~\cite{Tupciev1973} first  introduced the notion of self-similar viscosity approximations to the Riemann problem, in which the viscous profiles are sought to depend on the ratio $x/t$ of the space to the time coordinate. Next, the complete treatment of general hyperbolic systems by Dafermos~\cite{Dafermos1973Structure}, Dafermos and DiPerna~\cite{DafermosDiPerna1976}, and Tzavaras~\cite{Tzavaras1996} described also the wave structure of the limiting profiles. In particular, Tzavaras obtained a general framework, leading to uniform variation estimates for sufficiently close Riemann states of a strictly hyperbolic system and relates the internal structure of shocks to traveling waves. Schecter, Plohr, and Marchesin~\cite{SchecterPlohrMarchesin2004} developed a numerical continuation method for the boundary-value problem associated with the Dafermos viscosity. Slemrod~\cite{Slemrod1989} and Fan and Slemrod~\cite{FanSlemrod1993} studied self-similar artificial viscosity in phase dynamics and systems of mixed hyperbolic--elliptic type. Joseph and LeFloch addressed viscous and discrete approximations, self-similar diffusion with boundary conditions, and vanishing relaxation~\cite{JosephLeFloch1999,JosephLeFloch2002II,JosephLeFloch2002III}. Joseph and LeFloch~\cite{JosephLeFloch2007} also treated physical viscosity and capillarity in phase dynamics; this latter work motivates us in the present paper in which we partly generalize~\cite{JosephLeFloch2007} to relativistic fluids and to a space-time viscous operator.

\subsection{Large Riemann data and the Euler limit}
\paragraph{Main results.}
Namely, we consider the Riemann problem when  the viscosity coefficient in \eqref{eq:covariant-velocity-viscosity} is taken of the form\footnote{slightly changing the notation} 
 be $\eps\tau_+(t,x)$, where $\eps>0$ is a constant viscosity parameter and $\tau_+(t,x)\coloneqq\sqrt{\max\{t^2-x^2,0\}}$ for $t>0$. Inside the future light cone, the Lorentzian radius $\tau=\sqrt{t^2-x^2}$ replaces the time factor in Dafermos's self-similar approximation, and the coefficient vanishes outside the cone. The resulting model is given in Section~\ref{subsection:self-similar-velocity-viscosity} below. The construction of two half-profiles joined at a regular intermediate state follows the approach of Joseph and LeFloch~\cite{JosephLeFloch2007} for physical viscosity.

For genuinely nonlinear pressure laws satisfying the sound-speed bounds in \eqref{eq:convex-pressure-assumptions} below, we establish three properties.
\begin{itemize}
\item We construct global viscous profiles for arbitrary Riemann states, with positive lower density bounds and uniform bounds on density, rapidity, and total variation (cf.~Theorem~\ref{thm:convex-global-profiles} below).
\item We prove convergence of the full family to the classical Euler Riemann solution, with the prescribed initial trace and every smooth convex entropy inequality. For fixed data the differences converge in $L^1(d\zeta)$; for convergent data in a compact set they converge in this space after subtraction of the corresponding Riemann steps (cf.~Corollary~\ref{cor:compact-riemann-data} below).
\item We show that the viscous profiles equal the prescribed states outside $|x|\le t$. In the Euler limit the nonconstant region lies in a smaller cone $|x|\le t\tanh Z_0$, with finite $Z_0$. The bounds and viscosity threshold are uniform on compact sets of Riemann data.
\end{itemize}

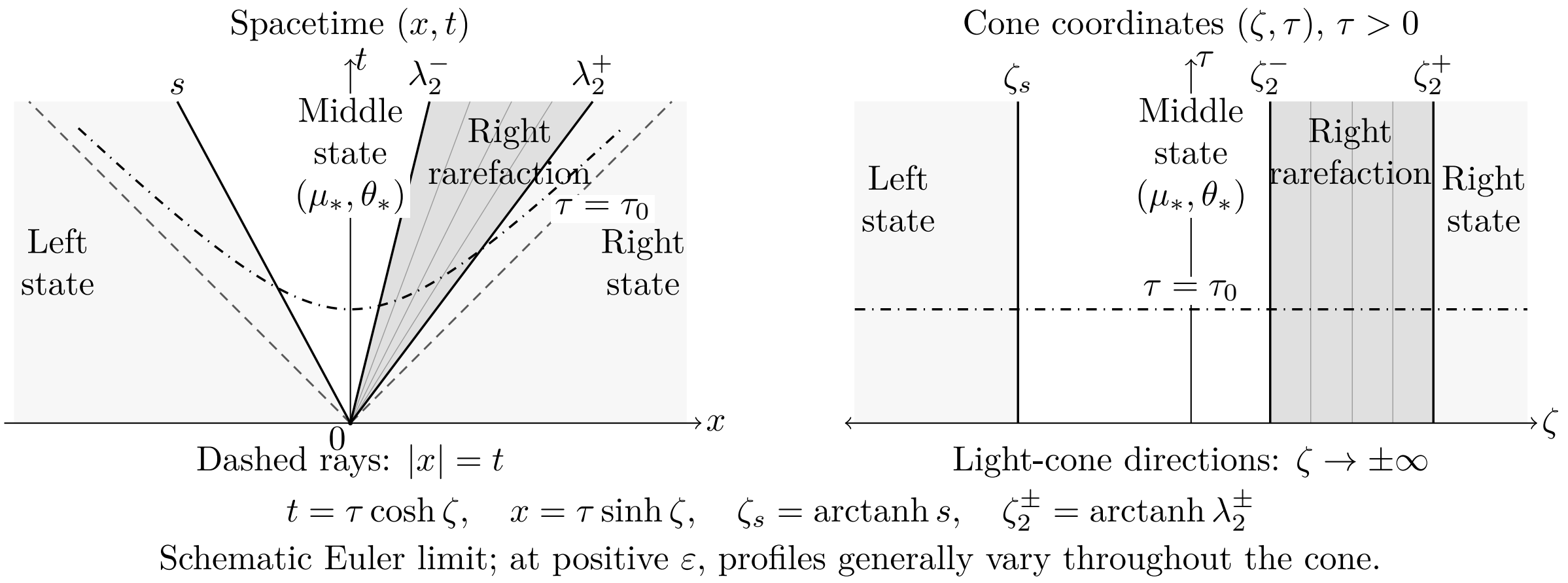
\begin{figure}[htbp]
\centering
\begingroup
\begin{tikzpicture}[
  x=1cm,y=1cm,
  every node/.style={font=\fontsize{10}{12}\selectfont,inner sep=1.1pt},
  axis/.style={->,line width=.4pt},
  lightcone/.style={densely dashed,line width=.55pt,black!65},
  leftedge/.style={line width=.65pt,black},
  rightedge/.style={line width=.65pt,black},
  rightinside/.style={line width=.25pt,black!40},
  taulevel/.style={dash dot,line width=.65pt,black}
]
  \begin{scope}[shift={(3.55,0)}]
    \fill[black!3] (-3.4,0)--(0,0)--(-1.75,3.25)--(-3.4,3.25)--cycle;
    \fill[black!3] (0,0)--(3.4,0)--(3.4,3.25)--(2.45,3.25)--cycle;
    \fill[black!12] (0,0)--(.8,3.25)--(2.45,3.25)--cycle;
    \foreach \xx in {1.21,1.63,2.04}
      \draw[rightinside] (0,0)--(\xx,3.25);
    \draw[lightcone] (0,0)--(-3.25,3.25);
    \draw[lightcone] (0,0)--(3.25,3.25);
    \draw[leftedge] (0,0)--(-1.75,3.25);
    \draw[rightedge] (0,0)--(.8,3.25);
    \draw[rightedge] (0,0)--(2.45,3.25);
    \draw[taulevel,domain=-2.75:2.75,samples=100]
      plot (\x,{sqrt(1.15*1.15+\x*\x)});
    \node[anchor=north west,fill=white,text=black]
      at (2.02,2.31) {$\tau=\tau_0$};
    \draw[axis] (-3.5,0)--(3.55,0) node[right] {$x$};
    \draw[axis] (0,0)--(0,3.68) node[right] {$t$};
    \fill (0,0) circle[radius=.65pt];
    \node[below left] at (0,0) {$0$};
    \node[align=center] at (-2.96,1.63) {Left\\state};
    \node[align=center] at (2.96,1.63) {Right\\state};
    \node[align=center] at (1.61,2.75) {Right\\rarefaction};
    \node[fill=white,align=center] at (0,2.70)
      {Middle\\state\\$(\mu_*,\theta_*)$};
    \node[above] at (-1.75,3.28) {$s$};
    \node[above] at (.8,3.28) {$\lambda_2^-$};
    \node[above] at (2.45,3.28) {$\lambda_2^+$};
    \node at (0,4.02) {Spacetime $(x,t)$};
    \node at (0,-.39) {Dashed rays: $|x|=t$};
  \end{scope}
  \begin{scope}[shift={(12.05,0)}]
    \fill[black!3] (-3.4,0) rectangle (-1.75,3.25);
    \fill[black!3] (2.45,0) rectangle (3.4,3.25);
    \fill[black!12] (.8,0) rectangle (2.45,3.25);
    \foreach \xx in {1.21,1.63,2.04}
      \draw[rightinside] (\xx,0)--(\xx,3.25);
    \draw[leftedge] (-1.75,0)--(-1.75,3.25);
    \draw[rightedge] (.8,0)--(.8,3.25);
    \draw[rightedge] (2.45,0)--(2.45,3.25);
    \draw[taulevel] (-3.4,1.15)--(3.4,1.15);
    \draw[<->,line width=.4pt] (-3.5,0)--(3.5,0) node[right] {$\zeta$};
    \draw[axis] (0,0)--(0,3.68) node[right] {$\tau$};
    \node[fill=white,text=black] at (0,1.34) {$\tau=\tau_0$};
    \node[align=center] at (-2.96,2.27) {Left\\state};
    \node[align=center] at (2.96,2.27) {Right\\state};
    \node[align=center] at (1.61,2.75) {Right\\rarefaction};
    \node[fill=white,align=center] at (0,2.70)
      {Middle\\state\\$(\mu_*,\theta_*)$};
    \node[above] at (-1.75,3.28) {$\zeta_s$};
    \node[above] at (.8,3.28) {$\zeta_2^-$};
    \node[above] at (2.45,3.28) {$\zeta_2^+$};
    \node at (0,4.02) {Cone coordinates $(\zeta,\tau)$, $\tau>0$};
    \node at (0,-.39) {Light-cone directions: $\zeta\to\pm\infty$};
  \end{scope}
  \node at (7.8,-.89)
    {$t=\tau\cosh\zeta,\quad x=\tau\sinh\zeta,\quad
      \zeta_s=\operatorname{arctanh}s,\quad\zeta_2^\pm=\operatorname{arctanh}\lambda_2^\pm$};
  \node at (7.8,-1.39)
    {Schematic Euler limit; at positive $\varepsilon$, profiles generally vary throughout the cone.};
\end{tikzpicture}
\endgroup
\caption{A first-family shock and a second-family rarefaction, separated by a constant intermediate state, shown in spacetime and in cone coordinates. The shock and the two edges of the rarefaction are rays in $(x,t)$ and vertical lines in $(\zeta,\tau)$. The level $\tau=\tau_0$ is a hyperbola in spacetime; the dashed rays are $|x|=t$.}
\label{fig:overview-cone-coordinates}
\end{figure}

\begin{figure}[htbp]
\centering
\begingroup%
\begin{tikzpicture}[x=1cm,y=1cm,
  font=\fontfamily{lmr}\fontsize{10}{12}\selectfont,
  every node/.style={inner sep=1pt},
  line cap=round,line join=round]
\path[use as bounding box] (0,0) rectangle (16.6,4.7);
\begin{scope}
  \node at (4.05,4.35) {(a) Density (schematic)};
  \fill[black!5] (2.4,1) rectangle (4.65,3.85);
  \draw[->,thin] (.55,1) -- (7.8,1) node[anchor=west] {$\zeta$};
  \draw[->,thin] (.55,.95) -- (.55,3.9) node[anchor=south] {$\mu$};
  \node[anchor=east] at (.48,1) {$0$};
  \foreach \xx in {2.4,4.65,6.25}
    \draw[black!35,densely dashed] (\xx,.92) -- (\xx,3.85);
  \draw[black,line width=1pt]
    (.85,1.65) -- (2.4,1.65) -- (2.4,2.3) -- (4.65,2.3)
    .. controls (5.15,2.3) and (5.7,3.2) .. (6.25,3.2) -- (7.4,3.2);
  \node at (1.15,1.95) {$\mu_l$};
  \node at (6.95,3.5) {$\mu_r$};
  \node at (3.525,2.6) {$\mu_*>0$};
  \node[align=center] at (3.525,3.35) {Middle\\state};
  \node[anchor=north] at (2.4,.84) {$\zeta_s$};
  \node[anchor=north] at (4.65,.84) {$\zeta_2^-$};
  \node[anchor=north] at (6.25,.84) {$\zeta_2^+$};
  \node[anchor=north] at (2.1,.38) {left shock};
  \node[anchor=north] at (5.8,.38) {right rarefaction};
\end{scope}
\begin{scope}[xshift=8.35cm]
  \node at (4.05,4.35) {(b) Rapidity (schematic)};
  \fill[black!5] (2.4,1) rectangle (4.65,3.85);
  \draw[->,thin] (.55,1) -- (7.8,1) node[anchor=west] {$\zeta$};
  \draw[->,thin] (.55,.95) -- (.55,3.9) node[anchor=south] {$\theta$};
  \foreach \xx in {2.4,4.65,6.25}
    \draw[black!35,densely dashed] (\xx,.92) -- (\xx,3.85);
  \draw[black,line width=1pt]
    (.85,2.45) -- (2.4,2.45) -- (2.4,1.6) -- (4.65,1.6)
    .. controls (5.15,1.6) and (5.7,3.4) .. (6.25,3.4)
    -- (7.4,3.4);
  \node at (1.15,2.75) {$\theta_l$};
  \node at (6.95,3.68) {$\theta_r$};
  \node at (3.525,1.95) {$\theta_*$};
  \node[anchor=north] at (2.4,.84) {$\zeta_s$};
  \node[anchor=north] at (4.65,.84) {$\zeta_2^-$};
  \node[anchor=north] at (6.25,.84) {$\zeta_2^+$};
  \node[anchor=north] at (2.1,.38) {left shock};
  \node[anchor=north] at (5.8,.38) {right rarefaction};
\end{scope}
\end{tikzpicture}%
\endgroup%
\caption{Schematic Euler density and rapidity profiles for a first-family shock at $\zeta_s$ followed by a second-family rarefaction on $[\zeta_2^-,\zeta_2^+]$. Across the shock, density increases and rapidity decreases; both increase across the rarefaction. The intermediate state $(\mu_*,\theta_*)$ is constant and has $\mu_*>0$. At positive $\eps$, viscous profiles generally vary throughout the light cone.}
\label{fig:overview-profiles}
\end{figure}
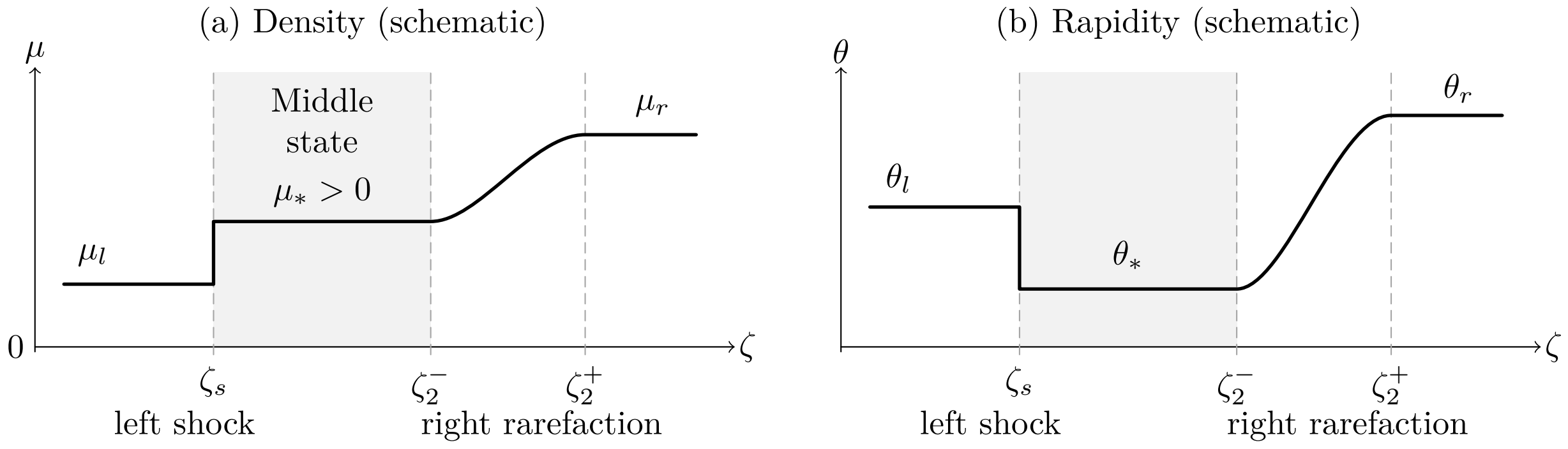

\paragraph{Outline of this paper.}
Section~\ref{section==2} derives the profile equations and constructs local profiles. Section~\ref{section:euler-away-vacuum} proves global existence and the vanishing-viscosity limit. Section~\ref{section:traveling-waves} identifies the traveling waves governing isolated compression layers, and Section~\ref{section==7} presents the numerical experiments. Figures~\ref{fig:overview-cone-coordinates} and~\ref{fig:overview-profiles} show a solution containing a first-family shock wave and by a second-family rarefaction wave.

\section{The viscous Riemann problem}
\label{section==2}

\subsection{The Lorentz-invariant model in hyperboloidal coordinates}
\label{subsection:self-similar-velocity-viscosity}

\paragraph{The divergence model.}
We define a Lorentz-invariant viscosity using the proper time from the Riemann origin. The resulting divergence form preserves constant states and admits self-similar profiles.
Inside the future light cone from the origin, namely $t>|x|$, we define the Lorentzian radius
\bse
\be
\label{eq:cone-radius}
\tau\coloneqq\sqrt{t^2-x^2}.
\ee
For $t>0$, set $\tau_+(t,x)\coloneqq\sqrt{\max\{t^2-x^2,0\}}$. We choose the conservative cone-scaled equation as our spacetime model:
\be
\label{eq:conservative-tau-model-main}
\del_\alpha T^{\alpha\beta}(\mu,U)
=
\eps\,\del_\alpha
\bigl(\tau_+\del^\alpha U^\beta\bigr).
\ee
\ese
The weight $\tau_+$ is homogeneous of degree one under spacetime dilations, so both sides have the same scaling, and every constant state solves this equation. Inside the cone, the product rule gives

\bse
\be
\label{eq:weighted-divergence-product}
\del_\alpha\bigl(\tau\del^\alpha U^\beta\bigr)
=\tau\Box U^\beta+\del_\alpha\tau\,\del^\alpha U^\beta.
\ee
The cone rapidity $\zeta\coloneqq\operatorname{arctanh}(x/t)$ is orthogonal to $\tau$ in the Lorentzian metric, since
\be
\label{eq:cone-radius-rapidity-orthogonality}
\del_\alpha\tau\,\del^\alpha\zeta=0.
\ee
\ese
For a profile $U=U(\zeta)$, the second term in the product rule therefore vanishes, and the spacetime model reduces to
\be
\label{eq:tau-box-model-main}
\del_\alpha T^{\alpha\beta}(\mu,U)
=\eps\,\tau\Box U^\beta.
\ee
The placement of the weight matters: since $\Box\tau=-1/\tau$, the product rule gives
\be
\Box(\tau U^\beta)=\tau\Box U^\beta+2\del_\alpha\tau\,\del^\alpha U^\beta-\frac{U^\beta}{\tau}.
\ee
Thus $\Box(\tau U)=(U''-U)/\tau$ on a profile, whereas the weighted divergence equals $U''/\tau$. In particular, the weighted divergence preserves constant states, while $\Box(\tau U)=-U/\tau$ for constant $U$. The divergence form also gives the viscous flux used to extend the profiles across the light cone.

The weighted divergence is Lorentz covariant relative to the origin of the Riemann data: the stress tensor and the velocity derivatives transform tensorially, while $\tau_+$ is a Lorentz scalar. The factor $\tau$ is the light-cone analogue of the time factor in Dafermos-type self-similar approximations. The viscous system has a hyperbolic principal part. In the first-order formulation used in Section~\ref{section==7} below, its characteristic speeds are the two light speeds and the fluid velocity; thus they remain in $[-1,1]$ for smooth physical states.

For the Riemann solution, compact support of the perturbation is established by extending the interior profile by its endpoint states and checking that the normal flux has matching traces on the two light rays. The decay of the interior profile as $|\zeta|\to\infty$ is used in this check, since $\zeta=\operatorname{arctanh}(x/t)$ tends to infinity at a light ray. A spatial Laplacian instead produces parabolic diffusion. Already for the heat equation, a nontrivial step has a Gaussian tail outside every bounded interval at positive times.

\paragraph{The light-cone frame.}
Inside the future light cone, we use the Milne coordinates $(\tau,\zeta)$, given by
\bse
\be
t=\tau\cosh\zeta,
\qquad
x=\tau\sinh\zeta.
\ee
In null coordinates, the same change of variables reads
\be
t+x=\tau e^\zeta,
\qquad
t-x=\tau e^{-\zeta}.
\ee
\ese
In these coordinates, the weighted divergence acts on each inertial component of the velocity by
\be
\label{eq:weighted-divergence-cone-coordinates}
\del_\alpha\bigl(\tau\del^\alpha U^\beta\bigr)
=\frac1\tau\left(-\del_\tau\bigl(\tau^2\del_\tau U^\beta\bigr)
+\del_\zeta^2 U^\beta\right).
\ee
The variable \(\zeta\) is the rapidity of the spacetime ray.  The associated
orthonormal frame is
\bse
\be
R(\zeta)\coloneqq(\cosh\zeta,\sinh\zeta),
\qquad
S(\zeta)\coloneqq(\sinh\zeta,\cosh\zeta),
\ee
and its vectors satisfy
\(\ang{R}{R}=-1,\ \ang{S}{S}=1, \ \ang{R}{S}=0\).

We parametrize the fluid velocity by its rapidity:
\be
U(\theta)=(\cosh\theta,\sinh\theta),
\qquad
N(\theta)\coloneqq\del_\theta U(\theta)
=
(\sinh\theta,\cosh\theta).
\ee
The velocity and its normal then satisfy
\(\ang{U}{U}=-1, \ \ang{U}{N}=0, \ \ang{N}{N}=1\).
With the relative rapidity \(\delta\coloneqq\theta-\zeta\), the velocity
and normal vectors decompose in the light-cone frame as
\be
\label{eq:U-N-R-S-decomposition}
U=\cosh\delta\,R+\sinh\delta\,S,
\qquad
N=\sinh\delta\,R+\cosh\delta\,S.
\ee
\ese

\paragraph{Reduction of the wave operator.}
For a scalar profile \(f=f(\zeta)\), the wave operator satisfies
\bse
\be
\Box f=\frac1{\tau^2}f''.
\ee
For the vector \(U(\theta(\zeta))\), componentwise application of \(\Box\)
gives
\be
\label{eq:tau-box-U-main}
\tau\Box U
=
\frac1{\tau}
\left(
N\,\theta''
+
U\,(\theta')^2
\right).
\ee
\ese
The viscous profile equations are degenerate because only the rapidity occurs at second order, and the viscous term is nonlinear because it contains the square of the first derivative \(\theta'\).

\subsection{The profile equations and main theorem}

\paragraph{The coupled profile equations.}
A self-similar solution has density $\mu=\mu(\zeta)$ and rapidity $\theta=\theta(\zeta)$.
For the divergence-form model, the identity $\del_\alpha\tau\,\del^\alpha U=0$ holds throughout the interior cone for a self-similar profile. Multiplication by $\tau$ therefore gives the vector profile equation
\be
\label{eq:divergence-vector-profile}
S_\alpha(\zeta)\frac{d}{d\zeta}T^{\alpha\beta}(\mu,U)=\eps\frac{d^2U^\beta}{d\zeta^2},\qquad S_\alpha=(-\sinh\zeta,\cosh\zeta).
\ee

\begin{proposition}[Riemann profile equations]
The reduced Riemann profile equations associated with
\eqref{eq:tau-box-model-main} are
\bse
\be
\label{eq:rapidity-normal-main}
\mu'\sinh(\theta-\zeta)
+
h(\mu)\theta'\cosh(\theta-\zeta)
=
\eps(\theta')^2,
\ee
\be
\label{eq:rapidity-tangential-main}
c^2(\mu)\mu'\cosh(\theta-\zeta)
+
h(\mu)\theta'\sinh(\theta-\zeta)
=
\eps\theta''.
\ee
\ese
\end{proposition}

\begin{proof}
We write $\delta=\theta-\zeta$. Using the components $T^{00}=h\cosh^2\theta-p$, $T^{01}=h\cosh\theta\sinh\theta$, and $T^{11}=h\sinh^2\theta+p$, we project the left-hand side of \eqref{eq:divergence-vector-profile} on $U$ and $N$ to obtain
\be
\begin{aligned}
\bigl\langle U,S_\alpha (T^{\alpha\cdot})'\bigr\rangle
 &= -\mu'\sinh\delta-h\theta'\cosh\delta,\\
\bigl\langle N,S_\alpha (T^{\alpha\cdot})'\bigr\rangle
 &= c^2\mu'\cosh\delta+h\theta'\sinh\delta .
\end{aligned}
\ee
The corresponding projections of $U''$ in \eqref{eq:tau-box-U-main} are $-(\theta')^2$ and $\theta''$, since $\langle U,U\rangle=-1$, $\langle N,N\rangle=1$, and $\langle U,N\rangle=0$. This gives both profile equations.
\end{proof}

\paragraph{Pressure assumptions and uniform estimates.}
We assume that
\be
\label{eq:convex-pressure-assumptions}
\begin{gathered}
 p\in C([0,+\infty))\cap C^3((0,+\infty)),\qquad p(0)=0,\\
 0<b_-\le p'(\mu)\le b_+<1,\qquad \mathcal G(\mu)\coloneqq h(\mu)p''(\mu)+2p'(\mu)(1-p'(\mu))>0
 \quad(\mu>0),
 \end{gathered}
\ee
where $b_-,b_+$ are fixed constants, so that $\sqrt{b_-}\le c(\mu)\le\sqrt{b_+}<1$. The last condition is genuine nonlinearity of the two characteristic families. The state estimates use only the sound-speed bounds. Genuine nonlinearity gives a uniform upper bound for the rapidity derivative on the resulting compact state set, which need not be smaller than one, and determines the admissibility of the limiting shocks.

\begin{theorem}[The viscous Riemann problem for the relativistic Euler equations]
\label{thm:convex-global-profiles}
Consider the relativistic Euler equations~\eqref{eq:relativistic-euler} in Minkowski spacetime and their viscous form~\eqref{eq:conservative-tau-model-main}. Suppose that the pressure satisfies~\eqref{eq:convex-pressure-assumptions}. Fix any two Riemann states with densities $\mu_l,\mu_r>0$ and rapidities $\theta_l,\theta_r\in\mathbb R$. There are constants $\eps_0,m,M,\Theta,C>0$, depending only on these states and the pressure, and profiles $(\mu^\eps,\theta^\eps)$ for $0<\eps<\eps_0$, with the following properties.
\begin{enumerate}
\item Each profile solves the viscous profile equations on $\mathbb R$, has the prescribed limits at $\zeta=\pm\infty$, and satisfies $\mu^\eps\in C^1(\mathbb R)$ and $\theta^\eps\in C^2(\mathbb R)$.
The two half-profiles meet at a regular junction $\zeta_*^\eps$, where $\theta^\eps(\zeta_*^\eps)=\zeta_*^\eps$ and $(\theta^\eps)'(\zeta_*^\eps)=0$. Its extension by the constant state $z_l$ for $x<-t$ and $z_r$ for $x>t$ solves the viscous system in divergence form in the sense of distributions on $t>0$ and has the prescribed Riemann initial trace. For each $t>0$, its difference from the initial Riemann step is supported in $[-t,t]$. At fixed $\eps>0$, the approach to the constant states from inside the light cone is superexponential in $\zeta$. Thus there are $a_\eps,C_\eps>0$ such that
\be
\label{eq:light-cone-superexponential}
 |\mu^\eps(\zeta)-\mu_{l,r}|+|\theta^\eps(\zeta)-\theta_{l,r}|
 \le C_\eps\exp\bigl(-a_\eps e^{|\zeta|}\bigr)
\ee
for sufficiently large $|\zeta|$, with the left state used as $\zeta\to-\infty$ and the right state as $\zeta\to+\infty$.

\item The profiles satisfy the uniform estimates
\bse
\be
\label{eq:convex-uniform-state-bounds}
 0<m\le\mu^\eps\le M,\qquad |\theta^\eps|\le\Theta,
\ee
\be
\label{eq:convex-uniform-BV-energy}
 \TV(\mu^\eps)+\TV(\theta^\eps)\le C,\qquad
 \eps\int_{\mathbb R}|(\theta^\eps)'|^2\dd\zeta\le C.
\ee
\ese
\item For the prescribed Riemann data, the limit $(\mu,\theta)$ is the same for every choice of the constructed profiles. As $\eps\downarrow0$, one has
\be
\label{eq:convex-profile-convergence}
 \mu^\eps-\mu\longrightarrow0,\qquad
 \theta^\eps-\theta\longrightarrow0
 \quad\hbox{strongly in } L^1(\mathbb R).
\ee
The limit defines a weak solution with bounded variation of the relativistic Euler system, with the prescribed Riemann initial data. Its non-constant support lies inside the light cone: for each $t>0$, the support of the spatial distributional derivative of the limiting state is contained in $|x|\le t\tanh Z_0$ for some finite $Z_0$ independent of $t$. The limit satisfies all convex entropy inequalities stated in Proposition~\ref{prop:convex-entropy-limit}.
\end{enumerate}
\end{theorem}

\begin{itemize}
\item \textbf{Physical Riemann data.} The positive lower density bound prevents vacuum in the constructed profiles. The prescribed velocities $\tanh\theta_l$ and $\tanh\theta_r$ may be arbitrarily close to the speed of light, with constants and a range of $\eps$ that depend on the data.

\item \textbf{Uniformity of the estimates.} The constants in the theorem depend on the prescribed states and the pressure, through the global bounds $b_-,b_+$ and the resulting compact density interval. Corollary~\ref{cor:compact-riemann-data} below gives uniform constants and a uniform threshold on compact sets of Riemann data. The continuation argument also uses bounds on the rapidity derivative at fixed $\eps$, while the passage to the Euler limit uses the uniform state, total variation, and energy estimates of the theorem.

\item \textbf{Isothermal pressure.} For $p(\mu)=k^2\mu$ with $0<k<1$, the enthalpy is $h(\mu)=(1+k^2)\mu$ and the sound speed is $c=k$. The Riemann invariants are $R=\theta+\frac{k}{1+k^2}\log\mu$ and $S=\theta-\frac{k}{1+k^2}\log\mu$, and the characteristic levels satisfy $\theta-\zeta=\pm\operatorname{arctanh}k$.
\end{itemize}

\paragraph{The first-order equations.}
We construct solutions of \eqref{eq:rapidity-normal-main}--\eqref{eq:rapidity-tangential-main} by joining two half-profiles. The first equation recovers the density; the rapidity derivative satisfies a scalar equation with a quadratic term.
For the elastodynamic system with physical viscosity and zero capillarity studied by Joseph and LeFloch~\cite[Section~2]{JosephLeFloch2007}, the self-similar equations and their scalar reduction are
\be
-yv'-\sigma(w)'=\eps v'',\qquad -yw'-v'=0,
 \qquad
 \eps(yw')'+(y^2-\sigma_w(w))w'=0.
\ee
Once $w$ is known, the constraint $v'=-yw'$ determines $v$ up to an additive constant. In the relativistic model, the normal equation determines the density derivative away from the junction, and its quadratic term contributes to the equation for $q$.

We write $q\coloneqq d\theta/d\zeta$ for the derivative of the fluid rapidity and recall that $\delta=\theta-\zeta$. The matrix multiplying $(\mu',\theta')$ on the left of the profile equations has determinant $h(\mu)(\sinh^2\delta-c^2(\mu)\cosh^2\delta)$. Its zeros occur at the characteristic levels $\tanh\delta=\pm c(\mu)$. At the junction $\delta=0$, the fluid velocity agrees with the velocity of the spacetime ray, and the coefficient of $\mu'$ in the normal equation vanishes.

Away from the junction, the density $\mu$ and rapidity derivative $q=\theta'$ satisfy the non-autonomous first-order system
\bse
\be
\label{eq:mu-prime-rapidity-main}
\aligned
 \mu' & =\frac{\eps q^2-h(\mu)q\cosh\delta}{\sinh\delta},
\qquad\quad
  q'  =\frac1\eps A(\mu,\theta,\zeta)q+B(\mu,\theta,\zeta)q^2,
\endaligned
\ee
where
\be
\label{eq:A-coefficient-main}
 A\coloneqq h(\mu)\frac{\sinh^2\delta-c^2(\mu)\cosh^2\delta}{\sinh\delta},
 \qquad B\coloneqq c^2(\mu)\coth\delta.
\ee
\ese
Along a profile, we choose a positive integrating factor $G$ with $G'/G=-Bq$. Then $(Gq)'=AGq/\eps$, and hence, for a constant $C$,
\be
\label{eq:theta-prime-kernel-main}
 q(\zeta)=\frac{C}{G(\zeta)}
   \exp\left(\frac1\eps\int_{\zeta_0}^{\zeta}A(\mu(s),\theta(s),s)\,\dd s\right).
\ee
The zeros of $A$ locate the characteristic levels. The additional factor $G^{-1}$ comes from the normalization of the velocity; its behavior must be controlled near the junction.

\paragraph{Compatibility at the junction.}
This first-order system is singular at $\delta=0$, so we return to the original profile equations to identify the admissible junction values before constructing the two branches.
At the junction, the original coupled equations give
\be
\label{eq:junction-compatibility-main}
 q_*\bigl(\eps q_*-h(\mu_*)\bigr)=0.
\ee
We select the regular junction $q_*=0$ and construct its two sides with the prescribed limiting states.

\subsection{Linearized profiles and branch monotonicity}
\label{subsection:linearized-fluid-profiles}
\label{subsection:nonlinear-monotonicity}

\paragraph{Linearized profiles.}
The linearized profiles provide the matching map used to start the nonlinear construction. Their constant extensions outside the light cone also illustrate the finite propagation property of the viscous system.
We fix a constant state $(\mu_0,\theta_0)$ and put $h_0=h(\mu_0)$ and $c_0=c(\mu_0)\in(0,1)$. For perturbations $m,v$ and $\delta=\theta_0-\zeta$, the linearized equations are
\be
\label{eq:linearized-fluid-system}
 m'\sinh\delta+h_0v'\cosh\delta=0,\qquad
 c_0^2m'\cosh\delta+h_0v'\sinh\delta=\eps v''.
\ee
The quadratic term has zero linearization at a constant state, so this linear system determines how changes in the junction values and amplitudes affect the prescribed endpoint states.

\begin{proposition}[Linearized profiles and their matching map]
\label{thm:linearized-fluid-profiles}
For $0<\eps<h_0c_0^2$ and arbitrary endpoint perturbations $(m_l,v_l)$, $(m_r,v_r)$, the system~\eqref{eq:linearized-fluid-system} has a unique solution with these limits and $m\in C^1$, $v\in C^2$. Its matching matrix is invertible and remains uniformly invertible as $\eps\to0$.
\end{proposition}
\begin{proof}
Eliminating $m'$ gives $\eps(v')'=A_0v'$, with $A_0=h_0((1-c_0^2)\sinh\delta-c_0^2/\sinh\delta)$. We define
\bse
\be
 H_\eps(\delta)\coloneqq
 e^{-h_0(1-c_0^2)\cosh\delta/\eps}
 |\tanh(\delta/2)|^{h_0c_0^2/\eps},\qquad
 Z_\eps\coloneqq\int_0^\infty H_\eps(s)\,ds,
\ee
\be
 M_\eps\coloneqq Z_\eps^{-1}\int_0^\infty\coth s\,H_\eps(s)\,ds>1.
\ee
\ese
The function $H_\eps(\theta_0-\zeta)$ solves the scalar equation for $v'$. We normalize it by $Z_\eps$ separately on the two half-lines to obtain kernels $\varphi_-,\varphi_+$ of integral one. Every solution has $v'=\alpha_-\varphi_-+\alpha_+\varphi_+$ and $m'=-h_0\coth\delta\,v'$. Therefore
\be
 \begin{pmatrix}v_r-v_l\\m_r-m_l\end{pmatrix}
 =\begin{pmatrix}1&1\\-h_0M_\eps&h_0M_\eps\end{pmatrix}
   \begin{pmatrix}\alpha_-\\\alpha_+\end{pmatrix}.
\ee
The determinant is $2h_0M_\eps>0$. The kernels decay faster than exponentially at infinity and behave at zero as a smooth positive factor times $|\delta|^\beta$, $\beta=h_0c_0^2/\eps>1$. Thus the integrals are finite, the recovered solution has the stated regularity, and the half-line representation proves uniqueness. We set $\phi(r)\coloneqq\eps\log H_\eps(r)$ for $r>0$. Its derivative is $h_0(c_0^2\cosh^2r-\sinh^2r)/\sinh r$, so the phase has a unique maximizer $r_c=\operatorname{arctanh}c_0$, with $\phi''(r_c)<0$. To apply Laplace's method to the numerator of $M_\eps$, choose $d\in(0,r_c)$ so small that $r\phi'(r)\ge\kappa>0$ for $0<r\le d$. Then $H_\eps(r)\le e^{\phi(d)/\eps}(r/d)^{\kappa/\eps}$. Since $\coth r\le C/r$ on this interval, the contribution from $(0,d)$ is bounded by $C(\eps/\kappa)e^{\phi(d)/\eps}$. This is exponentially smaller than the contribution at $r_c$. On $[d,+\infty)$, the weight $\coth r$ is bounded and the phase tends to minus infinity exponentially in $r$, so Laplace's method gives $M_\eps\to\coth r_c=1/c_0$. This proves uniform invertibility at the constant state.
\end{proof}

\paragraph{Monotonicity on each branch.}
Monotonicity of the density and rapidity on each side of the junction controls their total variation. The state estimates below make this bound uniform.

\begin{proposition}[Monotonicity of the density on a regular branch]
\label{prop:nonlinear-monotonicity}
Consider a smooth solution of the profile system, with $h>0$,
on the two sides of a regular junction $\zeta_*$, where
$\delta(\zeta_*)=0$ and $q(\zeta_*)=0$.
Assume that $\delta>0$ on the left interval and $\delta<0$ on the
right interval. On each interval the rapidity is either constant or
strictly monotone, and the density is monotone. If $q>0$, then
\be
\eps q<h(\mu)\cosh\delta.
\ee
If the state stays in a fixed compact state subset, the total variation of the state on the union of the two branches is bounded independently
of $\eps$ by twice the sum of the density and rapidity ranges of that set.
\end{proposition}

\begin{proof}
Away from the junction the equation
$q'=Aq/\eps+Bq^2$ implies that $q$ either vanishes identically or has a fixed sign on each connected branch, by uniqueness for the
scalar linear equation $q'=(A/\eps+Bq)q$ along the solution. If $q<0$, the reconstruction identity gives
\be
\mu'=\frac{q(\eps q-h\cosh\delta)}{\sinh\delta},
\ee
whose sign is fixed on either branch.

If $q>0$, define $L\coloneqq\eps q/(h\cosh\delta)>0$.
At a point where $L=1$, the reconstruction identity gives $\mu'=0$.
The second profile equation then gives $q'/q=h\sinh\delta/\eps$.
Since $\delta'=q-1$, differentiation of $L$ at such a point yields
\be
\frac{L'}{L}
 =\frac{q'}q-\frac{h_\mu\mu'}h-\tanh\delta\,(q-1)
 =\tanh\delta.
\ee
At the regular junction, $L$ tends to zero. On the left interval,
any excursion above $1$ would require a downward crossing of $1$
before reaching the junction, whereas the derivative at every crossing is positive. On the right interval, an upward crossing
from the initial value zero is excluded because the derivative at $L=1$ is negative. Thus $L<1$ on both branches. The reconstruction
identity again makes the sign of $\mu'$ fixed, proving monotonicity. If $q\equiv0$, the first equation gives $\mu'\equiv0$ away from the
junction. The total-variation estimate follows by adding the endpoint increments of the two monotone branches.
\end{proof}

\paragraph{The particle current and the derivative estimate.}
We fix $\bar\mu>0$ and define
\bse
\be
n(\mu)\coloneqq\exp\left(\int_{\bar\mu}^{\mu}\frac{dm}{h(m)}\right).
\ee
The identity $n'=n/h$ and the first profile equation give
\be
\label{eq:profile-current-identity}
\bigl(n(\mu)\sinh(\theta-\zeta)\bigr)'
 +n(\mu)\cosh(\theta-\zeta)
 =\eps\frac{n(\mu)}{h(\mu)}(\theta')^2.
\ee
\ese
Thus profiles in a fixed compact subset satisfy $\eps\int_I|\theta'|^2\,d\zeta\le C_I$ on each compact interval $I$.

\subsection{Construction and matching at the junction}
\label{subsection:fixed-epsilon-construction}

\paragraph{The regular junction.}
We seek solutions of the viscous equations~\eqref{eq:conservative-tau-model-main} in the sense of distributions on $t>0$, with the prescribed Riemann initial trace. At an interior ray $\zeta=\zeta_*$, two smooth half-profiles with continuous velocity solve the equations across the ray if and only if the normal fluxes agree:
\be
\label{eq:viscous-junction-flux}
\bigl[S_\alpha T^{\alpha\beta}-\eps(U^\beta)'\bigr]_{\zeta_*}=0 .
\ee
The brackets denote the right trace minus the left trace. Indeed, the vector profile equation has the form $(S_\alpha T^{\alpha\beta}-\eps(U^\beta)')'=S'_\alpha T^{\alpha\beta}$. Its right-hand side is locally integrable and the only possible atom on the left is the jump of the flux. The prime denotes differentiation in $\zeta$.

It remains to analyze the point
\(\delta=\theta-\zeta=0\).
The original equations have smooth coefficients there, but the coefficient of $\mu'$ in the first one vanishes. The construction in this section selects the regular alternative in~\eqref{eq:junction-compatibility-main},
\be
q(\zeta_*)=0,
\qquad
\theta(\zeta_*)=\zeta_*.
\ee
At this regular junction the common velocity is tangent to the ray, so $S_\alpha U^\alpha=0$. The trace of $(U^\beta)'$ is zero on each side. Consequently \eqref{eq:viscous-junction-flux} reduces to $[p(\mu)]_{\zeta_*}=0$, which is equivalent to equality of the two density traces since $p'>0$. Thus the construction uses a common density $\mu_*$ and a common rapidity $\theta_*=\zeta_*$ together with zero rapidity derivative on both sides.

\begin{figure}[htbp]
\centering
\begingroup
\begin{tikzpicture}[
  x=1cm,y=1cm,
  every node/.style={font=\fontsize{10}{12}\selectfont,inner sep=1.1pt},
  axis/.style={->,line width=.4pt},
  diagonal/.style={densely dashed,line width=.55pt,black!65},
  guide/.style={densely dotted,line width=.4pt,black!40},
  leftcurve/.style={line width=.85pt,blue!60!black},
  rightcurve/.style={line width=.85pt,orange!70!black},
  outward/.style={->,line width=.55pt,black!70}
]
  \begin{scope}[shift={(3.6,0)}]
    \draw[axis] (-2.85,-1.55)--(2.98,-1.55) node[right] {$\zeta$};
    \draw[axis] (-2.85,-1.55)--(-2.85,1.5) node[above] {$\theta$};
    \draw[guide] (-2.85,0)--(0,0)--(0,-1.55);
    \draw (-2.9,0)--(-2.8,0);
    \draw (0,-1.6)--(0,-1.5);
    \node[left] at (-2.91,0) {$\theta_*$};
    \node[below] at (0,-1.6) {$\zeta_*$};
    \draw[diagonal] (-1.35,-1.35)--(1.39,1.39);
    \node[above] at (1.22,1.43) {$\theta=\zeta$};
    \draw[guide] (-2.65,-.95)--(-1.95,-.95);
    \draw[guide] (1.95,.95)--(2.65,.95);
    \draw[leftcurve,<-,domain=-2.65:0,samples=100]
      plot (\x,{.95*tanh((\x/1.8)^3)});
    \draw[rightcurve,->,domain=0:2.65,samples=100]
      plot (\x,{.95*tanh((\x/1.8)^3)});
    \draw[line width=.4pt,black!65] (-.44,0)--(.44,0);
    \fill (0,0) circle[radius=1.2pt];
    \draw[line width=.35pt,black!55] (.08,-.07)--(.83,-.65);
    \node[align=center,fill=white] at (1.42,-.93)
      {$(\zeta_*,\theta_*)$\\$\theta_*=\zeta_*$};
    \node[text=blue!60!black] at (-1.73,-.34) {$\delta>0$};
    \node[text=orange!70!black] at (1.87,.49) {$\delta<0$};
    \node at (-2.35,-.73) {$\theta_l$};
    \node at (2.35,1.15) {$\theta_r$};
    \node at (0,1.99) {Regular junction};
    \draw[outward] (-.22,-2.08)--(-2.42,-2.08)
      node[midway,above] {$r$ increases};
    \draw[outward] (.22,-2.08)--(2.42,-2.08)
      node[midway,above] {$r$ increases};
    \node at (0,-2.49) {$q_*=\theta'(\zeta_*)=0$};
  \end{scope}
\end{tikzpicture}
\endgroup
\caption{The regular junction $(\zeta_*,\theta_*)$ of the two half-profiles satisfies $\theta_*=\zeta_*$ and $q_*=\theta'(\zeta_*)=0$. The diagonal is dashed and the two halves have different colors. The drawing illustrates an increasing profile schematically.}
\label{fig:overview-junction}
\end{figure}
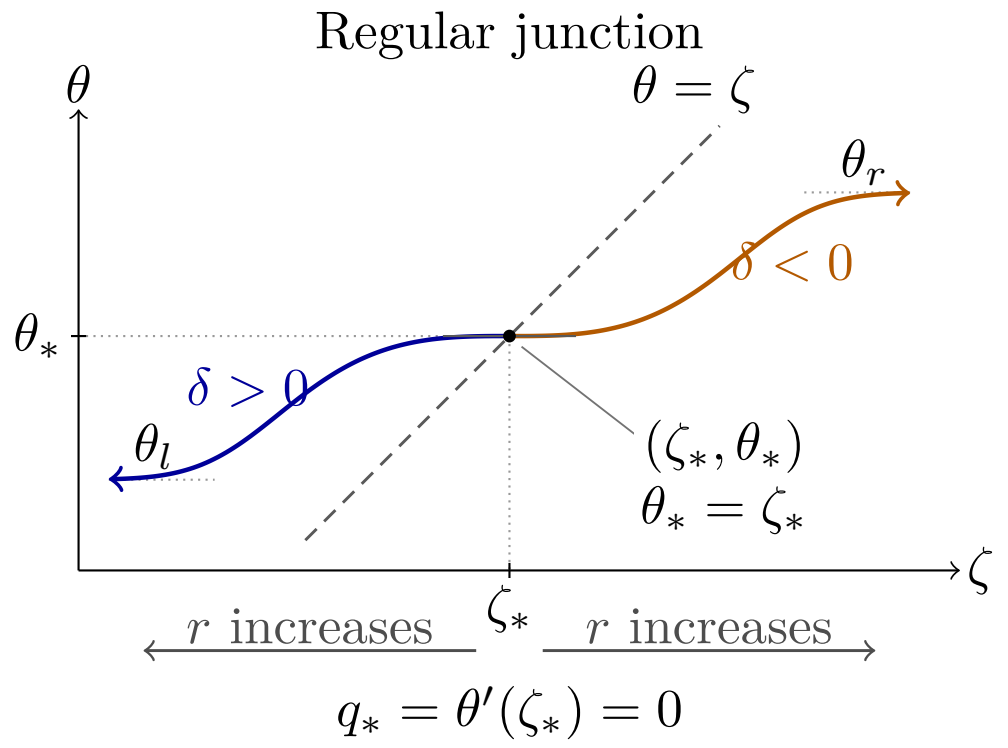

Figure~\ref{fig:overview-junction} shows the intersection of the rapidity profile with the diagonal $\theta=\zeta$. The two sides of this point are constructed separately, with $r=|\theta-\zeta|$ increasing away from the junction.

\paragraph{The equations in the distance to the junction.}
The equations in $r$ have the same form on both sides of the junction. This lets us solve the singular initial-value problem for both branches with the same local construction.

\begin{lemma}[Regular form at the junction]
\label{lemma:regular-junction-complete}
Let \(I\) be a punctured neighborhood of a point \(\zeta_*\) where
\(\delta(\zeta_*)=0\).  Assume that \(|q|\le1/2\) on \(I\).  Then \(\delta\)
is a valid local coordinate, and the profile system can be written as
\bse
\be
\label{eq:mu-delta-complete}
\frac{d\mu}{d\delta}
=
\frac{q}{1-q}
\frac{h(\mu)\cosh\delta-\eps q}{\sinh\delta},
\ee
\be
\label{eq:q-delta-complete}
\delta \frac{dq}{d\delta}
=
\Acal(\mu,\delta,q)q,
\ee
where the coefficient in the second equation is
\be
\label{eq:Acal-junction-complete}
\Acal(\mu,\delta,q)
\coloneqq
\frac{1}{\eps(1-q)}
\left[
h(\mu)
\frac{\delta\bigl(c^2(\mu)\cosh^2\delta-\sinh^2\delta\bigr)}{\sinh\delta}
-
\eps c^2(\mu)q\,\delta\coth\delta
\right].
\ee
At the regular point \((\delta,q)=(0,0)\), this coefficient satisfies
\be
\Acal(\mu,0,0)
=
\frac{h(\mu)c^2(\mu)}{\eps}>0.
\ee
\ese
\end{lemma}

\begin{proof}
Since \(\delta'=q-1\) and \(|q|\le1/2\), the variable \(\delta\) is a local coordinate, and the normal equation~\eqref{eq:rapidity-normal-main} gives
\bse
\be
\mu'
=
\frac{\eps q^2-h(\mu)q\cosh\delta}{\sinh\delta}.
\ee
Dividing by \(\delta'=q-1\) gives \eqref{eq:mu-delta-complete}.  Next,
\eqref{eq:rapidity-tangential-main} gives
\be
q'
=
\frac{1}{\eps}
\bigl(
c^2(\mu)\mu'\cosh\delta
+
h(\mu)q\sinh\delta
\bigr).
\ee
\ese
Substituting the expression for \(\mu'\) and dividing again by
\(\delta'=q-1\) gives \eqref{eq:q-delta-complete}.  The value at
\((\delta,q)=(0,0)\) follows from
\(\frac{\delta}{\sinh\delta}\to1, \qquad \delta\coth\delta\to1\).\qedhere
\end{proof}

\paragraph{Local existence and prescribed traces.}
We construct the two half-profiles at a regular junction and match their prescribed limiting states. The local result below also gives the starting point for the global continuation argument.

\begin{proposition}[Local nonlinear profiles for fixed $\eps$]
\label{thm:fixed-epsilon-local}
Let $p\in C^3$ near $\mu_0>0$, with $h>0$ and $0<p'<1$ there,
and fix $\theta_0\in\mathbb R$. For each fixed \(0<\eps<h(\mu_0)p'(\mu_0)\),
there is a neighborhood $\mathcal U_\eps$ of $z_0=(\mu_0,\theta_0)$ such that every pair
$z_l,z_r\in\mathcal U_\eps$ is connected by a profile on $\mathbb R$ with a regular junction. The profile has
$\mu\in C^1$, $\theta\in C^2$, takes its values in a fixed  neighborhood of $z_0$, and satisfies
\be
\TV(\mu,\theta)\le C_\eps|z_r-z_l|.
\ee
The construction gives a locally unique profile in the junction
parameter family described below. The neighborhood $\mathcal U_\eps$
and the constant $C_\eps$ may depend on $\eps$.
\end{proposition}

\begin{proof}
We construct each half-profile in the coordinate
$r=|\delta|$, starting at $r=0$, with junction density $\mu_*$
near $\mu_0$ and junction position $\zeta_* =\theta_*$ near
$\theta_0$. Throughout the proof $\eps$ is fixed. The equations for
$\mu$ and $q=\theta'$ are the same on both half-lines:
\bse
\be
\label{eq:fuchsian-fixed-epsilon}
\frac{d\mu}{dr}=\frac qr f(\mu,r,q),
 \qquad
 r\frac{dq}{dr}=\Acal(\mu,r,q)q,
\ee
where the first coefficient is
\be
f(\mu,r,q)\coloneqq
 \frac{r(h(\mu)\cosh r-\eps q)}{(1-q)\sinh r},
\ee
and the second coefficient is
\be
\Acal(\mu,r,q)\coloneqq
 \frac1{\eps(1-q)}
 \left(
 h(\mu)\frac{r(c^2(\mu)\cosh^2r-\sinh^2r)}{\sinh r}
 -\eps c^2(\mu)q\,r\coth r
 \right).
\ee
\ese
The quotients extend smoothly in $r$ at zero. We set
$\beta_*\coloneqq h(\mu_*)p'(\mu_*)/\eps$ and shrink the junction
density neighborhood so that $\beta_*>1$ throughout it.
The value $q=0$, $\mu=\mu_*$ is the constant half-profile.

\begin{enumerate}
\item \textbf{Integral equations at the junction.}
For a small real amplitude $C$, set $q(r)=r^{\beta_*}v(r)$.
The required half-profile is determined near zero by
\be
\label{eq:fuchsian-integral-map}
\begin{aligned}
 \mu(r)&=\mu_*+
 \int_0^r s^{\beta_*-1}v(s)
 f(\mu(s),s,s^{\beta_*}v(s))\,ds,\\
 v(r)&=C\exp\left(
 \int_0^r\frac{\Acal(\mu(s),s,s^{\beta_*}v(s))-\beta_*}{s}\,ds
 \right).
 \end{aligned}
\ee
We give the local estimates needed for this integral construction.
On a fixed compact state neighborhood and for $|q|<1/2$, the coefficient
satisfies
\bse
\be
|\Acal(\mu,r,q)-\beta_*|
 \le C_\eps\bigl(|\mu-\mu_*|+r^2+|q|\bigr).
\ee
For continuous bounded $v$, the first integral equation is a
contraction in the norm
$\sup_{0<r\le d}r^{-\beta_*}|\mu(r)-\mu_*|$
when $d$ is small and $v$ remains in a fixed small ball. Indeed its
Lipschitz constant with respect to $\mu$ is at most
$C_\eps d^{\beta_*}\|v\|_\infty$.
It yields
$|\mu(r)-\mu_*|\le C_\eps\|v\|_\infty r^{\beta_*}$
and a Lipschitz dependence on $v$ in this weighted norm.
After this equation is solved, the second equation maps a small
closed ball of continuous functions $v$ into itself. Its exponent
is bounded by
\(C_\eps\bigl(d^2+d^{\beta_*}\|v\|_\infty\bigr)\),
and the Lipschitz constant of the map is bounded by
$C_\eps|C|d^{\beta_*}$ on that ball. We choose $d$ and $|C|$ so
these contraction constants are less than one. This proves existence
and uniqueness of the local half-profile with
\be
q(r)=Cr^{\beta_*}(1+o(1)),
 \qquad \mu(r)-\mu_*=\cO(r^{\beta_*}).
\ee
\ese
\item \textbf{Dependence on the junction parameters.} We justify the $C^1$ dependence on $(\mu_*,C)$ at the singular endpoint. To work in a space independent of the junction density, we write $\mu(r)=\mu_*+r^{\beta_*}M(r)$. After the change of integration variable $s=rt$, the first integral equation becomes
\be
M(r)=\int_0^1 t^{\beta_*-1}v(rt) f\bigl(\mu_*+(rt)^{\beta_*}M(rt),rt,(rt)^{\beta_*}v(rt)\bigr)\,dt.
\ee
This map acts on a fixed ball of $C([0,d])$ and has the same contraction bound as the weighted equation. Its parameter derivatives are continuous there: differentiating the powers introduces only $\log t$ under the integral and factors of the form $(rt)^{\beta_*}\log(rt)$, extended by zero at $rt=0$. The weight $t^{\beta_*-1}(1+|\log t|)$ has an integrable majorant independent of the parameters when $\beta_*$ ranges over a compact interval in $(1,+\infty)$. Thus $M$ depends continuously differentiably on $(\mu_*,v)$. Differentiating the reconstruction at fixed $v$ gives
\bse
\be
\del_{\mu_*}\mu(r)-1
 =\cO_\eps\bigl(r^{\beta_*}(1+|\log r|)\bigr).
\ee
The coefficient $\Acal$ is continuously differentiable in $(\mu,q)$, with a locally Lipschitz density derivative. Since $\mu-\mu_*=\cO(r^{\beta_*})$, differentiation at fixed $v$ gives
\be
\left|\del_{\mu_*}
 \bigl(\Acal(\mu(r),r,r^{\beta_*}v(r))-\beta_*\bigr)\right|
 \le C_\eps\left(
 r^{\beta_*}(1+|\log r|)+r^2\right).
\ee
\ese
After division by $r$, all terms are integrable at zero. The derivatives with respect to $v$ have the integrable bound $C_\eps r^{\beta_*-1}$ in the second exponent. After solving the first integral equation for $\mu$, the second map acts on the fixed space $C([0,d])$ of functions $v$. On a smaller parameter neighborhood, $\beta_*$ is bounded below by a constant greater than one. The preceding bounds justify differentiation under the integrals and continuity of the derivatives with an integrable majorant independent of the parameters. The contraction theorem therefore gives $C^1$ dependence of the fixed points.

\item \textbf{Extension to the half-line.}
We next extend the half-profile to $r=+\infty$ for small $C$.
We choose a compact density interval around $\mu_0$ on which
$h\ge h_*>0$ and $c\le c^*<1$. On each fixed interval $[d,R]$, ordinary differential-equation dependence on initial conditions keeps $\mu$ close to $\mu_*$ and $q$ small when $C$ is small. We choose $R$ sufficiently large that, throughout this compact density
interval and for $|q|\le1/2$, the coefficient satisfies
\bse
\be
\frac{\Acal(\mu,r,q)}r\le-b_\eps e^r
 \quad(r\ge R)
\ee
for some $b_\eps>0$. This follows from $c^*<1$ in the explicit
formula for $\Acal$. As long as the solution remains in the interval,
\be
|q(r)|\le |q(R)|
 \exp\bigl(-b_\eps(e^r-e^R)\bigr),
 \qquad |\mu_r|\le C_\eps|q(r)|.
\ee
\ese
The total density change on $[R,+\infty)$ is therefore
$\cO_\eps(|q(R)|)$. Reducing $|C|$ keeps the density away from the
boundary of the chosen interval. These estimates close the
continuation argument and prove existence of $\mu(+\infty)$.
To differentiate the endpoint values, write the half-profile equations as $\mu_r=F(\mu,r,q)$ and $q_r=G(\mu,r,q)$. On the same tail, the explicit coefficients give
\be
 |F_\mu|\le K|q|,\qquad |F_q|\le K,\qquad
 |G_\mu|\le Ke^r|q|,\qquad G_q\le-\lambda e^r
 \qquad(r\ge R),
\ee
for positive constants $K,\lambda$ at fixed $\eps$, after increasing $R$ and restricting $|q|$ if necessary. In particular, the leading term of $G_q$ is $-h(1-p')e^r/[2\eps(1-q)^2]$. These estimates hold throughout a neighborhood of the constant profile in the parameter space. Let $(\dot\mu,\dot q)$ denote a derivative with respect to either $\mu_*$ or $C$, which exists on each finite interval by the local construction and the ordinary differential equations. The variational equation for $q$ gives
\be
 |\dot q(r)|\le e^{-\lambda(e^r-e^R)}|\dot q(R)|
 +K\int_R^r e^{-\lambda(e^r-e^s)}e^s|q(s)|\,|\dot\mu(s)|\,ds.
\ee
Integrating up to any finite radius $r$ and using
$\int_s^\infty e^{-\lambda(e^t-e^s)}\dd t\le e^{-s}/\lambda$ gives
\be
 \int_R^r|\dot q(t)|\dd t
 \le K_1|\dot q(R)|+K_1\int_R^r |q(s)|\,|\dot\mu(s)|\dd s
\ee
for some $K_1>0$. The equation for $\dot\mu$ now gives a Gronwall bound with integrable coefficient proportional to $|q|$. Thus $\dot\mu$ is uniformly bounded. Substitution in the bound for $\dot q$, using the decay of $q$, gives
$|\dot q(r)|\le K_2\exp(-\lambda_1(e^r-e^R))$ for some $K_2>0$ and $0<\lambda_1<\min(\lambda,b_\eps)$. Indeed, after the substitution $u=e^s$ the convolution of the two exponentials is bounded by a linear factor in $e^r-e^R$ times an exponential, and reducing the rate absorbs that factor.

The derivatives of $\mu_r$ and $q/(1-q)$ therefore have an integrable majorant independent of the parameters on the tail. Near zero, the weighted junction construction gives an integrable bound by a multiple of $r^{\beta_0-1}(1+|\log r|)$ for some fixed $\beta_0>1$. Differentiation under the integrals and continuous dependence on finite intervals show that the limiting density and the rapidity increment are $C^1$ in $(\mu_*,C)$.

The rapidity coordinate and fluid rapidity are recovered by
\be
\begin{aligned}
 \zeta_-(r)&=\theta_*-\int_0^r\frac{ds}{1-q_-(s)},
 &\theta_-(r)&=\zeta_-(r)+r,\\
 \zeta_+(r)&=\theta_*+\int_0^r\frac{ds}{1-q_+(s)},
 &\theta_+(r)&=\zeta_+(r)-r.
 \end{aligned}
\ee
Here $q_-$ and $q_+$ use independent small amplitudes $C_-$ and $C_+$. Since $|q_\pm|<1/2$, both coordinate maps are invertible
onto their half-lines. The condition $\beta_*>1$ gives $C^1$ density and $C^2$ rapidity across the regular junction.

\item \textbf{Matching the far-field states.}
It remains to match the four prescribed limits. At the constant
profile the derivative of each $q$ with respect to its amplitude
is a positive multiple of $H_\eps(r)$ from Proposition~\ref{thm:linearized-fluid-profiles}. We may therefore use
\bse
\be
\alpha_\pm\coloneqq\int_0^{+\infty}
 \frac{q_\pm(r)}{1-q_\pm(r)}\,dr
\ee
as local amplitude coordinates: their derivatives with respect
to $C_\pm$ at zero are strictly positive and finite.
The derivative of the map assigning the limiting states to the junction parameters in the coordinates
$(\mu_*,\theta_*,\alpha_-,\alpha_+)$ is
\be
\begin{aligned}
 d\mu_l&=d\mu_*+h_0M_\eps\,d\alpha_-,
 &d\theta_l&=d\theta_*-d\alpha_-,\\
 d\mu_r&=d\mu_*+h_0M_\eps\,d\alpha_+,
 &d\theta_r&=d\theta_*+d\alpha_+.
 \end{aligned}
\ee
\ese
This map is invertible: its difference equations have the matching
matrix of Proposition~\ref{thm:linearized-fluid-profiles}, and the two
remaining equations determine the junction state. The inverse function theorem gives a locally unique solution for all boundary
states sufficiently close to $z_0$, with a neighborhood depending on $\eps$. The strengths vanish on the diagonal $z_l=z_r$ and
depend $C^1$ on the data, so they are bounded by
$C_\eps|z_r-z_l|$. The same bound holds for the density increments
by the half-profile estimates. Monotonicity on each branch then gives the total-variation bound in the proposition. This completes the proof of Proposition~\ref{thm:fixed-epsilon-local}.\qedhere
\end{enumerate}
\end{proof}

\section{Uniform estimates, global continuation, and the Euler limit}
\label{section:euler-away-vacuum}

\subsection{Uniform  bounds for matched half-profiles}

\paragraph{The continuation path.}
We specify the class of half-profiles used for global matching. State bounds and compactness of the matching parameters then allow continuation from constant data to the prescribed Riemann states.
We continue the local profiles to the prescribed Riemann states by combining the weighted flux identities with bounds on the rapidity derivative. A normalized derivative remains below one, while the rapidity derivative itself may cross one. We therefore use distance from the junction in cone rapidity for global continuation and use relative rapidity only near the junction and on the tails. The state and variation bounds then give compactness as $\eps$ tends to zero.
We use the pressure hypotheses \eqref{eq:convex-pressure-assumptions}. We continue to use
$h(\mu)=\mu+p(\mu)$ and $c^2(\mu)=p'(\mu)$, and put
$a_+\coloneqq1+b_+$. These assumptions imply
\be
\label{eq:convex-pressure-comparisons}
b_-\mu\le p(\mu)\le b_+\mu,\qquad
 (1+b_-)\mu\le h(\mu)\le a_+\mu.
\ee
The functions $h$ and $\mu-p(\mu)$ are strictly increasing.
The latter property uses $p'<1$; it does not require $p''\ge0$.

\begin{lemma}[The domain of the conservative variables]
\label{lemma:convex-conservative-domain}
Under \eqref{eq:convex-pressure-assumptions}, the map $(\mu,\theta)\mapsto W=(T^{00},T^{01})$ is a global diffeomorphism from $(0,+\infty)\times\mathbb R$ onto the convex cone
\be
\mathcal C\coloneqq\{(E,M)\in\mathbb R^2:E>|M|\}.
\ee

\end{lemma}

\begin{proof}
We write $u\coloneqq\tanh\theta$, so the conservative variables satisfy
\bse
\be
E=\frac{h}{1-u^2}-p,\qquad
 M=\frac{hu}{1-u^2},\qquad
 E\pm M=\frac{\mu\pm pu}{1\mp u}>0.
\ee
The last inequalities follow from $0<p<\mu$ and $|u|<1$.
Conversely, fix $(E,M)\in\mathcal C$. The scalar equation
\be
\label{eq:convex-conservative-inverse}
(E-\mu)(E+p(\mu))=M^2,\qquad 0<\mu\le E,
\ee
\ese
has exactly one solution. Its left-hand side decreases from $E^2$
to zero, since its derivative for $0<\mu<E$ is
$-(E+p)+(E-\mu)p'<0$. For $M=0$ the solution is $\mu=E$.
Setting $u=M/(E+p(\mu))$ then gives $|u|<1$ and recovers both
conservative equations. The Jacobian is $h(\cosh^2\theta-p'(\mu)\sinh^2\theta)>0$, which proves smooth invertibility with the regularity of the pressure.
\end{proof}

\paragraph{Matched half-profiles.}
We set $\ell\coloneqq|\zeta-\zeta_*|$ on each side of the regular junction, with $\zeta_*=\theta_*$. A dot denotes differentiation with respect to $\ell$. For either half-profile put
\bse
\be
v(\ell)\coloneqq\int_0^\ell q(s)\dd s,\qquad r(\ell)\coloneqq\ell-v(\ell).
\ee
The left and right states are recovered from $\zeta_\pm=\theta_*\pm\ell$ and $\theta_\pm=\theta_*\pm v$. On both sides $r=|\theta-\zeta|$, and the equations become
\be
\label{eq:global-branch-ell}
\dot\mu=\frac{q(h\cosh r-\eps q)}{\sinh r},\qquad
\dot r=1-q,\qquad
\dot q=\frac{q}{\eps\sinh r}\bigl(h(p'\cosh^2r-\sinh^2r)-\eps p'q\cosh r\bigr).
\ee
\ese
Unlike the equations in $r$, this system is regular at $q=1$ whenever $r>0$. We consider pairs which extend to $\ell=+\infty$, have $r>0$ there, positive finite limiting densities, $q\to0$, and finite increments
\be
\alpha_\pm\coloneqq\int_0^{+\infty}q_\pm(\ell)\dd\ell.
\ee
At the junction we require
\be
\label{eq:convex-junction-exponent}
\beta_*\coloneqq h(\mu_*)p'(\mu_*)/\eps>1.
\ee
The pair is matched when its endpoint states are the prescribed states, so that $d_l=\theta_*-\theta_l=\alpha_-$ and $d_r=\theta_r-\theta_*=\alpha_+$. Since $q$ has a fixed sign on each nonconstant half-profile, $v$ and the rapidity are monotone there, and $r=\ell-\alpha_\pm+o(1)$ at infinity.

\begin{lemma}[A normalized derivative and the inner part of each branch]
\label{lemma:normalized-derivative}
On a positive half-profile define $L\coloneqq\eps q/(h\cosh r)$ and $\kappa\coloneqq b_+/(1+b_+)<1$. Then $0<L<\kappa$. On both positive and negative branches the density is monotone, and its increment has the same sign as $q$.
If $r_b>0$ satisfies $\tanh^2r_b<b_-/(1+b_+)$, a positive matched branch passes through $0<r\le r_b$ only once, with $q<1$ and $\dot r>0$.
\end{lemma}
\begin{proof}\bse
With $b=p'(\mu)$, differentiation of $L$ using \eqref{eq:global-branch-ell} gives
\be
\label{eq:normalized-derivative-barrier}
\frac{\dot L}{L}
=\frac{h(1-L)}{\eps\sinh r}
 \bigl[\cosh^2r\,(b-(1+b)L)-\sinh^2r\bigr]-\tanh r.
\ee
At $L=\kappa$, one has $b-(1+b)\kappa=(b-b_+)/(1+b_+)\le0$, so $\dot L<0$ for every $r>0$. Since $L\to0$ at the junction, a first-crossing argument on each regular interval proves the bound. In particular $h\cosh r-\eps q>h\cosh r/(1+b_+)>0$ on a positive branch; on a negative branch positivity is immediate. The density equation proves the monotonicity assertion. This argument uses $h_\mu=1+b$ and the sound-speed bounds, without differentiating $b$.
On the other hand, in the indicated inner strip the coefficient of $q$ in its equation satisfies
\be
\frac{\dot q}{q}
=\frac{h\cosh^2r}{\eps\sinh r}\bigl(b(1-L)-\tanh^2r\bigr)>0.
\ee
Consequently the region $0<r\le r_b$, $q\ge1$ is forward invariant as long as the solution remains regular: $\dot q>0$ on $q=1$ and $\dot r\le0$ on $r=r_b$. A matched branch cannot enter this region, because its finite increment gives $r(\ell)\to+\infty$. Hence $q<1$ in the strip. A return through $r=r_b$ would require $\dot r\le0$ and therefore $q\ge1$, which has just been excluded. The initial passage is thus the only one. Near the junction, division by $\dot r$ gives
\be
\label{eq:convex-half-density}
\mu_r=\frac{q}{1-q}\frac{h(\mu)\cosh r-\eps q}{\sinh r}. \qedhere
\ee
\ese
\end{proof}

\paragraph{Weighted flux identities and upper bounds.}
We integrate the weighted flux on each branch between the junction and its limiting state. The resulting identities give one-sided bounds on the junction density and rapidity, according to the sign of the rapidity increment.

\begin{lemma}[Weighted fluxes for a variable sound speed]
\label{lemma:convex-null-fluxes}
Every such matched profile satisfies
\bse
\be
\label{eq:convex-null-fluxes}
\begin{aligned}
 p(\mu_*)e^{2\theta_*}
 &=\tfrac12h(\mu_l)e^{2\theta_l}
 -\int_{-\infty}^{\theta_*}(\mu-p(\mu))e^{2\zeta}\,d\zeta
 -\eps\int_{-\infty}^{\theta_*}e^{\theta+\zeta}q\,d\zeta,\\
 p(\mu_*)e^{-2\theta_*}
 &=\tfrac12h(\mu_r)e^{-2\theta_r}
 -\int_{\theta_*}^{+\infty}(\mu-p(\mu))e^{-2\zeta}\,d\zeta
 -\eps\int_{\theta_*}^{+\infty}e^{-\theta-\zeta}q\,d\zeta.
 \end{aligned}
\ee
For either side, let $(m,d)$ denote $(\mu_l,d_l)$ or $(\mu_r,d_r)$.
The following inequalities hold:
\be
\label{eq:convex-side-bounds}
\begin{aligned}
 h(\mu_*)&\le h(m)e^{-2d}, &&d\ge0,\\
 h(\mu_*)&\ge h(m)e^{-2d}, &&d\le0,\\
 p(\mu_*)&\le\tfrac12h(m)e^{-2d}-\tfrac12(m-p(m))
                    +\eps e^{-d}(-d), &&d\le0.
 \end{aligned}
\ee
\ese
\end{lemma}

\begin{proof}
We set $S_\alpha=(-\sinh\zeta,\cosh\zeta)$, $R_\alpha=(-\cosh\zeta,\sinh\zeta)$, and $P^\beta\coloneqq S_\alpha T^{\alpha\beta}-\eps(U^\beta)'$. The vector profile equation gives $P'=R_\alpha T^{\alpha\cdot}$ and hence
\be
(e^\zeta(P^0+P^1))'=-(\mu-p)e^{2\zeta}-\eps e^{\theta+\zeta}q,
 \qquad
 (e^{-\zeta}(P^0-P^1))'=-(\mu-p)e^{-2\zeta}-\eps e^{-\theta-\zeta}q.
\ee
At the junction, $P=p(\mu_*)S(\theta_*)$, which gives the values $\pm p(\mu_*)e^{\pm2\theta_*}$ for the corresponding weighted components. At the two exterior ends we use the identity
\be
e^{\pm\zeta}(P^0\pm P^1)
=h e^{\pm(\theta+\zeta)}\sinh(\theta-\zeta)\pm p e^{\pm2\zeta}
\mp\eps e^{\pm(\theta+\zeta)}q.
\ee
For the upper sign at the left end only the enthalpy term has a nonzero limit, namely $\tfrac12h(\mu_l)e^{2\theta_l}$. For the lower sign at the right end that term tends to $-\tfrac12h(\mu_r)e^{-2\theta_r}$. The weighted integrands are integrable on their respective half-lines because the state and the rapidity derivative have finite limits there and the weights decay exponentially. Integration on the left from $-\infty$ to the junction and on the right from the junction to $+\infty$ gives \eqref{eq:convex-null-fluxes}.

Divide the appropriate flux identity by $e^{\pm2\theta_*}$.
For $d\ge0$, the normalized $(\mu-p)$ integral lies between
$(\mu_*-p(\mu_*))/2$ and $(m-p(m))/2$, while the normalized
$q$ integral lies between zero and $d$. Keeping its favorable
sign gives the first inequality in \eqref{eq:convex-side-bounds}.
For $d\le0$, the order of the density bounds reverses and the
$q$ integral is nonpositive, of absolute value at most
$e^{-d}(-d)$. On the left, for example, $\theta-\theta_*\le-d$
and $\zeta-\theta_*\le0$, which bounds its exponential weight by
$e^{-d}$; on the right the same argument uses the opposite signs.
Using respectively the upper and lower bounds on the $(\mu-p)$
integral gives the remaining two inequalities.
\end{proof}

\paragraph{Lower density bound on a positive branch.}

Now, we fix $\bar\eps>0$ and the prescribed states, and define
\be
\label{eq:convex-upper-state-constants}
\begin{gathered}
 D\coloneqq\theta_r-\theta_l,\qquad D_+\coloneqq\max(D,0),\\
 \mu_{\min}\coloneqq\min(\mu_l,\mu_r),\qquad
 \mu_{\max}\coloneqq\max(\mu_l,\mu_r),\\
 h_{\min}\coloneqq h(\mu_{\min}),\qquad
 h_{\max}\coloneqq h(\mu_{\max}),\\
 M_p\coloneqq\max\left\{\mu_{\max},
 \frac{h(\mu_l)e^{2|D|}-(\mu_l-p(\mu_l))
            +2\bar\eps e^{|D|}|D|}{2b_-}\right\},\\
 \Theta_p\coloneqq\max(|\theta_l|,|\theta_r|)
                  +\frac12\log\frac{h_{\max}}{h_{\min}}.
 \end{gathered}
\ee
Every matched profile in the above class, for $0<\eps\le\bar\eps$,
satisfies $\mu\le M_p$ and $|\theta|\le\Theta_p$.
If $\mu_*\ge\mu_{\max}$, both increments are nonpositive and
$|d_l|,|d_r|\le|D|$. The third inequality of
\eqref{eq:convex-side-bounds} on the left, together with
$p(\mu_*)\ge b_-\mu_*$, bounds the density peak by $M_p$.
If $\mu_*\le\mu_{\min}$, the density upper bound is already an
endpoint bound, and both increments are nonnegative. In these
two cases $\theta_*$ lies between the prescribed rapidities.
In the remaining case the whole density profile lies in the endpoint
density interval. The first two inequalities of
\eqref{eq:convex-side-bounds} give
\be
|d_l|,|d_r|\le\frac12\log\frac{h_{\max}}{h_{\min}},
\ee
which bounds $\theta_*$ by $\Theta_p$. Monotonicity on each branch
extends these bounds to the full profile.

Since $b_-(1-b_+)>0$, continuity allows us to choose $r_b>0$
and $\gamma>0$ such that $\tanh^2r_b<b_-/(1+b_+)$ and
\bse
\be
\label{eq:convex-inner-radius}
b_-(\cosh^2r-b_+\cosh r)-\sinh^2r\ge\gamma
 \quad(0\le r\le r_b).
\ee
We then choose $0<r_a<r_b$ and set
\be
\label{eq:convex-density-floor-constants}
\begin{gathered}
 c_p\coloneqq\gamma/b_+,\qquad
 \chi_p\coloneqq
 \left(\frac{\tanh(r_a/2)}{\tanh(r_b/2)}\right)^{c_p}<1,\\
 L_{0,p}\coloneqq
 \frac{a_+\cosh r_a\,\chi_p}{c_p(1-\chi_p)},\qquad
 L_{1,p}\coloneqq a_+\coth r_a.
 \end{gathered}
\ee
\ese

\begin{lemma}[Density bound on a positive half-profile]
\label{lemma:convex-density-floor}
Let a positive matched half-profile satisfy
\eqref{eq:convex-junction-exponent}, have limiting density $m>0$,
and have increment $d=\int_0^{+\infty}q(\ell)\dd\ell$. Then
\be
\label{eq:convex-half-density-floor}
\mu_*\ge m\exp(-L_{0,p}-L_{1,p}d).
\ee
\end{lemma}

\begin{proof}
Lemma~\ref{lemma:normalized-derivative} permits the use of $r$ up to the first exit from the inner strip, with $q<1$ there and no later return. Density monotonicity implies $h(\mu)\ge h(\mu_*)$, and the junction
condition gives $\eps/h(\mu)<p'(\mu_*)\le b_+$.
Writing $b=p'(\mu)$ only within this proof, the equation for the rapidity derivative on a half-profile becomes
\bse
\be
q_r=\frac{h q}{\eps(1-q)\sinh r}
       \left(b\cosh^2r-\sinh^2r-\frac{\eps}{h}bq\cosh r\right).
\ee
The parenthesis is bounded below by the left-hand side of
\eqref{eq:convex-inner-radius}, since $b\ge b_-$ and
$\cosh^2r-b_+\cosh r>0$. Moreover $h/\eps>1/b_+$.
These bounds imply
\be
q_r\ge\frac{c_p}{\sinh r}\frac{q}{1-q}
 \quad(0<r\le r_b).
\ee
Integration of $(\log q-q)_r\ge c_p/\sinh r$ up to $r_b$,
using $\log q(r_b)-q(r_b)\le-1$, yields
\be
q(r)\le
 \left(\frac{\tanh(r/2)}{\tanh(r_b/2)}\right)^{c_p},
 \qquad q(r)\le\chi_p\quad(0<r\le r_a).
\ee
The density equation and \eqref{eq:convex-pressure-comparisons} give
\be
(\log\mu)_r\le a_+\frac{q}{1-q}\coth r.
\ee
\ese
On $(0,r_a]$, we use $\coth r\le\cosh r_a/\sinh r$, $1-q\ge1-\chi_p$, and the power bound. Since $(\tanh(r/2)^{c_p})_r=c_p\tanh(r/2)^{c_p}/\sinh r$, integration gives exactly the upper bound $a_+\cosh r_a\,\chi_p/[c_p(1-\chi_p)]=L_{0,p}$. After the first crossing of $r_a$, one has $r\ge r_a$ and $\dot\mu/\mu\le a_+q\coth r_a$. Integration in $\ell$ bounds the remaining logarithmic density increment by $L_{1,p}d$, proving \eqref{eq:convex-half-density-floor}.
\end{proof}

\paragraph{Uniform bounds.}

The lower density estimate on an expanding branch complements the bounds from the weighted fluxes. We combine them with monotonicity on each branch to obtain the same bounds throughout the continuation.

\begin{proposition}[Uniform bounds on the density and rapidity]
\label{prop:convex-uniform-states}
Define
\bse
\be
\label{eq:convex-uniform-floor}
m_p\coloneqq\min\{\mu_{\min},
        \mu_l\exp(-L_{0,p}-L_{1,p}D_+)\}>0.
\ee
All matched profiles of the above class with $0<\eps\le\bar\eps$
satisfy
\be
\label{eq:convex-uniform-states}
m_p\le\mu\le M_p,\qquad |\theta|\le\Theta_p.
\ee
The same constants apply along the data path
\be
\label{eq:convex-data-path}
\mu_r(s)=(1-s)\mu_l+s\mu_r,\qquad
 \theta_r(s)=(1-s)\theta_l+s\theta_r,
 \qquad 0\le s\le1.
\ee
\ese
\end{proposition}

\begin{proof}
Only the lower bound remains to be proved. If $\mu_*\ge\mu_{\min}$,
monotonicity already bounds the entire profile below by $\mu_{\min}$. Otherwise the rapidity derivative is positive on both open half-lines and
$0<d_l,d_r\le D_+$. Lemma~\ref{lemma:convex-density-floor} on the left bounds $\mu_*$ below by the second entry in
\eqref{eq:convex-uniform-floor}. Along \eqref{eq:convex-data-path}
the endpoint density and rapidity intervals shrink, while $|D(s)|\le|D|$ and $D(s)_+\le D_+$. The constants in the upper
and lower estimates can therefore be kept fixed.
\end{proof}

\subsection{Global matching of the two branches}

\paragraph{Bounds on the rapidity derivative.}
Genuine nonlinearity bounds the positive derivative uniformly in viscosity on a compact state set. The negative derivative remains bounded at each fixed viscosity, which suffices for continuation.

\begin{lemma}[Upper bound under genuine nonlinearity]
\label{lemma:convex-rapidity-derivative-bound}
Suppose the matched profiles take values in $[m,M]\times[-\Theta,\Theta]$, with $0<\eps\le\bar\eps$. Then
\be
\label{eq:convex-rapidity-derivative-gap}
-Q_\eps\le q\le Q_+,\qquad
Q_+\coloneqq\max\left\{1,\max_{m\le\mu\le M}\frac{2p'(\mu)(1-p'(\mu))}{\mathcal G(\mu)}\right\}<+\infty,
\ee
where $Q_\eps$ is finite at fixed viscosity and $Q_+$ is independent of viscosity.
\end{lemma}

\begin{proof}
The uniform bounds and the two monotone branches bound the total variation of $\theta$. Choose $Z>\Theta+\operatorname{arctanh}\sqrt{b_+}+2$. On $[-Z,Z]$, the identities $P'=R_\alpha T^{\alpha\cdot}$ and $P(\theta_*)=p(\mu_*)N(\theta_*)$ bound $P$. The relation $\eps U'=S_\alpha T^{\alpha\cdot}-P$ bounds $\eps|q|$ there. The exterior integrating-factor estimate in Lemma~\ref{thm:global-conditional-limit} bounds $q$ by its average on an adjacent unit interval; it uses only the compact state set and $p'\le b_+<1$. This proves a finite absolute bound at each fixed viscosity.

On a positive branch the maximum $Q$ is attained, since $q$ vanishes at zero and infinity. Put $b=p'(\mu)$, $b_\mu=p''(\mu)$, and
\bse
\be
\mathcal B\coloneqq h(\sinh^2r-b\cosh^2r)+\eps bq\cosh r.
\ee
The equation in $\ell$ is $\dot q=-q\mathcal B/(\eps\sinh r)$. At a maximum, $\mathcal B=0$ and $\dot{\mathcal B}\ge0$. With $K=h\cosh r-\eps Q>0$, differentiation at fixed $(r,q)$ gives
\be
\label{eq:convex-B-density-derivative}
\del_\mu\mathcal B=-\frac{\eps bQ\cosh r}{h}(1+b)-b_\mu\cosh r\,K.
\ee
The other partial derivative is $\del_r\mathcal B=2h(1-b)\sinh r\cosh r+\eps bQ\sinh r$. Using $\dot\mu=QK/\sinh r$ and $\dot r=1-Q$ therefore gives
\be
\label{eq:convex-full-maximum-inequality}
\begin{aligned}
 (1-Q)\sinh^2r\,[2h(1-b)\cosh r+\eps bQ]
 \ge \frac{\eps b(1+b)Q^2\cosh r}{h}K+b_\mu Q\cosh r\,K^2.
\end{aligned}
\ee
No division by $1-Q$ has been made. If $Q\ge1$ and $b_\mu\ge0$, the left side is nonpositive and the right side is positive, a contradiction. If $Q\ge1$ and $b_\mu<0$, substitute $\sinh^2r=b\cosh r\,K/h$ and collect the terms to obtain
\be
h\cosh r\,[2b(1-b)-Q\mathcal G]
\ge\eps Q\bigl[(b+2b^2-hb_\mu)Q-b^2\bigr]>0.
\ee
\ese
Thus $Q<2b(1-b)/\mathcal G$ in this case, which proves the asserted bound.
\end{proof}

\paragraph{Continuation from constant data.}
We start with the constant profile and vary the right state along a straight path. The estimates above control the parameters of every matched profile along this path, so the continuation argument can be carried out on a bounded parameter domain.

\begin{proof}[Construction of the viscous profiles in Theorem~\ref{thm:convex-global-profiles}]

\begin{enumerate}
\item \textbf{The endpoint map.}
We fix
$0<\bar\eps\le h(\mu_l)p'(\mu_l)$, depending only on the prescribed data and the pressure, take the constants $m_p,M_p,\Theta_p$ of
Proposition~\ref{prop:convex-uniform-states}, and choose a threshold
\bse
\be
\label{eq:convex-epsilon-threshold}
0<\eps_0<\min\{\bar\eps,\ b_-h(m_p),\ h(m_p)/2\}.
\ee
For each $\eps\in(0,\eps_0)$, the integral equations
\eqref{eq:fuchsian-integral-map}, which were derived for general
pressure, construct a local half-profile for every junction density
with $h(\mu_*)p'(\mu_*)>\eps$ and every finite amplitude
$C\in\mathbb R$. Near any fixed parameter pair, choose $C_{\max}>0$ exceeding the absolute values of the amplitudes in a bounded neighborhood and use the fixed ball $\|v\|_\infty\le2C_{\max}$. The exponent is bounded by $C_\eps(d^2+C_{\max}d^{\beta_*})$ and the contraction constant by $C_\eps C_{\max}d^{\beta_*}$. One sufficiently small radius $d$ therefore gives a contraction throughout that neighborhood, including when the limiting amplitude is zero.
The contraction estimates in
Section~\ref{subsection:fixed-epsilon-construction} give local
uniqueness and $C^1$ parameter dependence, with
\(q(r)=Cr^{\beta_*}(1+o(1)),\qquad \mu(r)-\mu_*=\cO(r^{\beta_*})\).
We denote by $\mathcal D_{\eps,p}$ the set of triples
$(\mu_*,C_-,C_+)$ with $h(\mu_*)p'(\mu_*)>\eps$ for which
both half-profiles extend globally in $\ell$ with $r>0$, have positive densities with positive finite limits, satisfy $q\to0$, and have finite increments $\alpha_\pm=\int_0^{+\infty}q_\pm\dd\ell$. The variable $r$ is used to construct the local germ and is not assumed to be a global coordinate.
We define
\be
\label{eq:convex-endpoint-map}
\Phi_{\eps,p}(\mu_*,C_-,C_+)
 \coloneqq(\mu_-^\infty,\mu_+^\infty,\alpha_-+\alpha_+).
\ee
\ese
For the prescribed states, the three matching conditions are
\be
\mu_-^\infty=\mu_l,\qquad
\mu_+^\infty=\mu_r,\qquad
\alpha_-+\alpha_+=\theta_r-\theta_l.
\ee
The last condition ensures that a single junction rapidity recovers both prescribed endpoint rapidities. Figure~\ref{fig:riemann-shooting-map} shows how the three parameters enter these conditions.

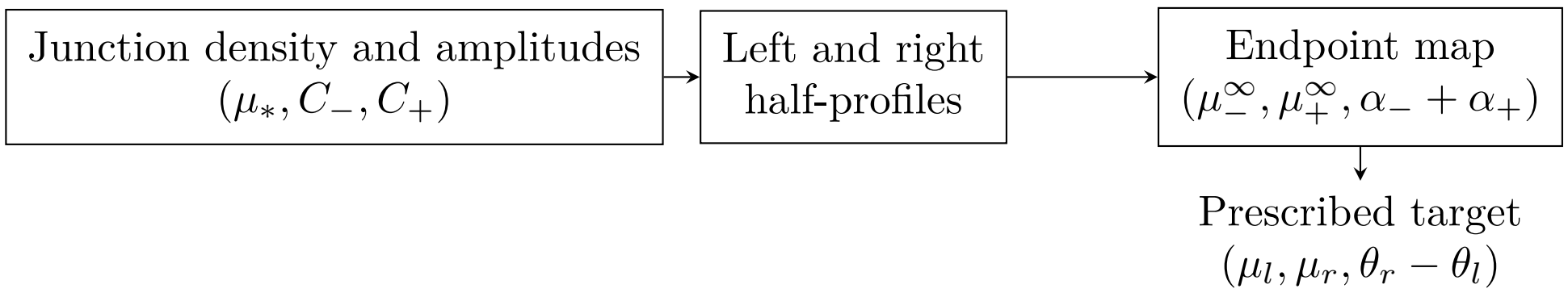
\begin{figure}[htbp]
\centering
\begin{tikzpicture}[every node/.style={font=\normalsize,align=center},>=stealth]
 \node[draw,inner sep=5pt] (params) at (0,0)
   {Junction density and amplitudes\\$(\mu_*,C_-,C_+)$};
 \node[draw,inner sep=5pt] (branches) at (4.25,0)
   {Left and right\\half-profiles};
 \node[draw,inner sep=5pt] (ends) at (8.4,0)
   {Endpoint map\\$(\mu_-^\infty,\mu_+^\infty,\alpha_-+\alpha_+)$};
 \draw[->] (params) -- (branches);
 \draw[->] (branches) -- (ends);
 \node (target) at (8.4,-1.35) {Prescribed target\\$(\mu_l,\mu_r,\theta_r-\theta_l)$};
 \draw[->] (ends) -- (target);
\end{tikzpicture}
\caption{The shooting parameters determine the two half-profiles. Matching means that the output of the endpoint map is the prescribed target.}
\label{fig:riemann-shooting-map}
\end{figure}

The target path is
$Y_s\coloneqq(\mu_l,\mu_r(s),\theta_r(s)-\theta_l)$, with
\eqref{eq:convex-data-path}. A preimage is matched by setting
$\theta_*\coloneqq\theta_l+\alpha_-$.

\textit{Continuity of the endpoint map.} Continuation on finite intervals uses \eqref{eq:global-branch-ell}, which is regular for $r>0$ even when $q=1$. On a matched branch the inner-strip lemma excludes a return to $r=0$. For an arbitrary admissible triple, positivity of $r$ on compact intervals and its divergence at infinity give the same local continuation margin. For an admissible triple, we take a compact positive density interval containing the limiting density in its interior. For $r\ge R$ and $|q|\le q_0<1$, the half-profile equations give $q_r/q\le-\lambda_\eps e^r$ and $|\mu_r|\le C_\eps|q|$ on this interval, where $\lambda_\eps>0$. The first bound follows from the leading term, proportional to $-h(1-p')e^r/\eps$. To reach this tail we continue in $\ell$ to the last crossing of a fixed large level $r=R$ of the given profile. There $|q|<1/2$ and $\dot r>1/2$, so the crossing time and its state depend continuously on the parameters. Nearby profiles enter the same tail intervals with a positive density margin. Only on this tail do we use $r$ again as an independent variable. The total remaining density drift is at most $C_\eps|q(R)|e^{-R}/\lambda_\eps$, which is smaller than that margin for large $R$. Thus the profiles remain in the intervals, and
\be
 |q(r)|\le |q(R)|\exp\bigl(-\lambda_\eps(e^r-e^R)\bigr)
 \qquad(r\ge R).
\ee
The admissible set is therefore open. The exponential bound dominates the integrands for the two endpoint quantities on the tail, while the junction estimates control their integrals near zero. Dominated convergence and continuous dependence on finite intervals give continuity of $\Phi_{\eps,p}$.

\item \textbf{Compactness of the matching parameters.}
\begin{itemize}
\item \textit{Bounds near the junction.} Every pre-image
of every $Y_s$ satisfies the uniform bounds of
Proposition~\ref{prop:convex-uniform-states}. Lemma
\ref{lemma:convex-rapidity-derivative-bound} gives, for this fixed $\eps$,
a uniform absolute bound on $q$. Increasing $Q_\eps$ if necessary, we use $|q|\le Q_\eps$ over all these preimages. On a positive branch Lemma~\ref{lemma:normalized-derivative} and the inner power estimate give $q\le\chi_p<1$ for $0<r\le r_a$; on a negative branch $q\le0$. Thus only the local calculation near the junction uses a separation from one.
At zero the coefficient in \eqref{eq:fuchsian-fixed-epsilon} is
\bse
\be
\Acal(\mu,0,q)
 =\frac{p'(\mu)(h(\mu)-\eps q)}{\eps(1-q)}.
\ee
It is bounded below by a positive constant on the compact local ranges of the state and $-Q_\eps\le q\le\chi_p<1$, since $h(m_p)>2\eps$.
There are therefore $0<r_0\le r_a$ and $\gamma_\eps>0$, uniform over the preimages, such that $\mathcal A\ge\gamma_\eps$ for $0<r\le r_0$. Each branch reaches this fixed level on its initial monotone part, and its first crossing occurs at $\ell\ge r_0/(1+Q_\eps)$ because $|\dot r|\le1+Q_\eps$. Thus the same initial $r$ interval is available for every preimage. The equation for $q$, followed by integration of $|\mu_r|\le C_\eps|q|/r$, gives
\be
|q(r)|\le Q_\eps(r/r_0)^{\gamma_\eps},\qquad
 |\mu(r)-\mu_*|\le C_\eps r^{\gamma_\eps}
 \quad(0<r\le r_0).
\ee
The coefficient $\Acal$ is even and smooth in $r$ at zero, so
\be
|\Acal(\mu(r),r,q(r))-\beta_*|
 \le C_\eps(r^{\gamma_\eps}+r^2).
\ee
Consequently the amplitude identity
\be
\label{eq:convex-amplitude-bound}
C_\pm=q_\pm(r_0)r_0^{-\beta_*}
 \exp\left(-\int_0^{r_0}
 \frac{\Acal(\mu_\pm(r),r,q_\pm(r))-\beta_*}{r}\,dr\right)
\ee
bounds both amplitudes uniformly, including the zero amplitude of a constant half-profile. Here
\be
1<\frac{b_-h(m_p)}{\eps}\le\beta_*
 \le\frac{b_+h(M_p)}{\eps},
\ee
\ese
and the integral in \eqref{eq:convex-amplitude-bound} is uniformly
bounded in absolute value.

\item \textit{Closure and convergence of the endpoint data.} Any sequence of preimages has a convergent parameter subsequence,
and its target parameters have a convergent subsequence in $[0,1]$.
The limiting exponent remains greater than one by
\eqref{eq:convex-epsilon-threshold}. The local integral construction uses one radius for the limit triple and all sufficiently late triples in the subsequence, since the amplitudes are bounded and the exponents range over a compact subset of $(1,+\infty)$. Continuous dependence at the end of this fixed interval gives initial values for ordinary continuation. On a sufficiently small fixed neighborhood of the junction the convergence is already given by the local construction. On a positive branch the estimate $q\le\chi_p$ before the first crossing of $r_a$, and the absence of a return afterward, give $r(\ell)\ge\min\{r_a,(1-\chi_p)\ell\}$. Negative branches have $r(\ell)\ge\ell$. The limiting profiles therefore stay uniformly away from $r=0$ on every interval separated from the junction. The limiting profiles inherit the uniform bounds and $-Q_\eps\le q\le Q_+$. The system \eqref{eq:global-branch-ell} is consequently regular on each remaining finite interval of $\ell$, even at $q=1$, and the profiles extend over all such intervals.

The increments along the target path are bounded by the state and rapidity bounds. Since $r=\ell-v$ and $|v|$ is bounded by the absolute increment of its monotone branch, $r\to+\infty$ uniformly as $\ell\to+\infty$. The explicit equation in $\ell$ gives
\bse
\be
\frac{\dot q}{q}=-\frac{h(\mu)(1-p'(\mu))}{2\eps}e^r+\cO_\eps(1),\qquad
|\dot\mu|\le C_\eps|q|
\ee
on a uniform tail, with the quotient denoting the coefficient also when $q=0$. Because $r=\ell+\cO(1)$ uniformly, integration gives
\be
|q(\ell)|\le C_\eps\exp\bigl(-\lambda_\eps(e^\ell-e^{\ell_0})\bigr),\qquad \ell\ge\ell_0,
\ee
\ese
for a fixed $\ell_0$ and positive $\lambda_\eps$. This bounds the integrands $q$ and $\dot\mu$ by an integrable function independent of the preimage. Passing first on finite intervals and then in the tails proves convergence of both limiting densities and increments. The limiting exponent remains greater than one, $r$ stays positive away from the junction, and the limiting densities remain at least $m_p$. The limit triple belongs to $\mathcal D_{\eps,p}$ and has the limiting target as its image. This proves compactness of the union of preimages inside the open admissible set.
\end{itemize}

\item \textbf{Continuation from the constant profile.}
\begin{itemize}
\item \textit{The constant profile and its linearization.} At $s=0$,
both prescribed densities equal $\mu_l$ and the total increment is
zero. If $\mu_*<\mu_l$, both increments are positive; if
$\mu_*>\mu_l$, both are negative. Both alternatives contradict
their zero sum. Equality forces both monotone half-profiles to
be constant. Thus the unique preimage is $(\mu_l,0,0)$.
The linearized construction of
Proposition~\ref{thm:linearized-fluid-profiles} applies to the frozen
values $h_l=h(\mu_l)$ and $c_l=\sqrt{p'(\mu_l)}$.
With its positive normalized moment $M_\eps$, the linearization, written in terms of the junction density and the first variations of the increments, has matrix
\be
D\Phi_{\eps,p}=
 \begin{pmatrix}
 1&h_lM_\eps&0\\
 1&0&h_lM_\eps\\
 0&1&1
 \end{pmatrix},\qquad
 \det D\Phi_{\eps,p}=-2h_lM_\eps\ne0.
\ee
To justify linearization of the endpoint map at the constant profile, consider a half-profile with junction density $\mu_*$ near $\mu_l$ and amplitude $C$. Proposition~\ref{thm:fixed-epsilon-local} gives $|q|<1/2$ globally for these small amplitudes. In this neighborhood $r$ is a global coordinate and $\int q\dd\ell=\int q/(1-q)\dd r$; the following calculation in $r$ is restricted to this neighborhood. On each finite interval, the local construction and continuous differentiability of the ODE flow give convergence of $q/C$ to the linearized half-profile as $(\mu_*,C)\to(\mu_l,0)$ with $C\ne0$. Near zero, $|q/C|\le C_\eps r^{\beta_0}$ for some $\beta_0>1$, by the junction integral equations. On the tail, the preceding estimate divided by $|C|$ gives a integrable bound independent of these parameters, since $q(R)=\cO(|C|)$ uniformly for $\mu_*$ near $\mu_l$. These bounds also dominate $q/((1-q)C)$ and $\mu_r/C$, the latter with an additional factor $1/r$ near zero. Dominated convergence thus gives the linearized increments and density changes, with errors $o(|C|)$ as $(\mu_*,C)\to(\mu_l,0)$. When $C=0$, the increment is zero and the limiting density equals $\mu_*$. This proves the asserted differentiability at the constant triple in the amplitude coordinates. The linearized coefficients taking $C_\pm$ to $\alpha_\pm$ are positive, so the derivative in amplitude coordinates is invertible as well.

\item \textit{Piecewise affine continuation.} We write $P$ for a parameter triple, $P_0$ for the constant triple, and $J\coloneqq D\Phi_{\eps,p}(P_0)$ for its invertible linearization. We choose a bounded polyhedral open neighborhood $\mathcal O\Subset\mathcal D_{\eps,p}$ of the compact union of preimages; it may be taken as the interior of a finite union of small cubes from a common grid. Set $F_s(P)\coloneqq\Phi_{\eps,p}(P)-Y_s$. These maps omit zero on $[0,1]\times\del\mathcal O$. Suppose $F_1$ has no zero in $\mathcal O$, so it also omits zero on $\overline{\mathcal O}$.

The expansion $F_0(P)=J(P-P_0)+o(|P-P_0|)$ permits a local interpolation to $J(P-P_0)$ using a cutoff equal to one near $P_0$ and supported in a small ball in $\mathcal O$. On that ball we have $|F_0(P)-J(P-P_0)|<|P-P_0|/(2\|J^{-1}\|)$ for $P\ne P_0$, so every interpolated map still has $P_0$ as its only zero. Concatenating the reverse of this interpolation with $F_s$ gives a continuous map $H$ on $[0,1]\times\overline{\mathcal O}$ whose final map is $F_1$. We choose $\rho>0$ so that its initial map is affine on $B_{2\rho}(P_0)$ and this closed ball lies in $\mathcal O$. The lateral boundary, the final face, and the initial face outside $B_\rho(P_0)$ have positive distances from zero under $H$. Choose $\eta>0$ smaller than one quarter of these distances and such that $\eta\|J^{-1}\|<\rho$.

We choose a compatible triangulation of $[0,1]\times\overline{\mathcal O}$ into finitely many simplices of dimension four, respecting its boundary faces, and interpolate $H$ affinely at their vertices. Uniform continuity makes the resulting map $\widehat H$ uniformly closer than $\eta$ to $H$ when the subdivision is fine enough. With mesh smaller than $\rho$, its initial map remains $J(P-P_0)$ on $B_\rho(P_0)$. For every $|\xi|<\eta$, the initial equation $\widehat H(0,P)=\xi$ consequently has the single solution $P=P_0+J^{-1}\xi$, while the lateral and final faces contain no solution.

\item \textit{Pairing the endpoints.} We choose such a $\xi$ outside the images of all faces of dimension at most two and of all simplices on which the affine map has rank below three. The excluded images lie in a finite union of affine sets of dimension at most two, which cannot cover a ball in $\mathbb R^3$. In a simplex of dimension four that meets $\widehat H^{-1}(\xi)$, the affine map has rank three, so its inverse image is a line intersected with the simplex. This intersection is a nondegenerate segment with endpoints in the relative interiors of faces of dimension three. Indeed, a singleton intersection would either lie in a face of dimension at most two or be tangent to a face of dimension three. The former case has been excluded; in the latter the line would lie in the affine span of that face, whose affine restriction would then have rank below three.

On every face of dimension three that is met, the affine restriction has rank three and hence the point mapping to $\xi$ is unique. Each interior face is shared by two simplices, each contributing one endpoint at this point, so the interior endpoints are paired. There are no endpoints on the lateral or final boundary, and there is precisely one on the initial face. This is impossible: a finite collection of segments with all interior endpoints paired has an even number of unpaired endpoints. The contradiction proves that $F_1$ has a zero in $\mathcal O$.
\end{itemize}

The recovery formulas $\zeta_\pm=\theta_*\pm\ell$ cover the entire rapidity line. Near the junction they agree with the local reconstruction in $r$, while finite increments give $r=\ell+\cO(1)$ on the tails.
The condition $\beta_*>1$ gives the claimed junction regularity.

\item \textbf{Uniform estimates and extension across the light cone.}
The two monotone
branches and \eqref{eq:convex-uniform-states} give, explicitly,
\be
\TV(\mu)\le2(M_p-m_p),\qquad
 \TV(\theta)\le4\Theta_p.
\ee
We choose a fixed interior rapidity interval containing all junctions
and the two unit strips used in the exterior estimates.
Integration of \eqref{eq:profile-current-identity}
bounds $\eps\int q^2$ on this interval: the boundary terms at the
endpoints and the integral of $n\cosh\delta$ are uniformly bounded,
and $n/h$ has a positive minimum on $[m_p,M_p]$.
The uniform exterior estimate~\eqref{eq:uniform-exterior-q}, proved in Section~\ref{subsection:global-conditional-limit} below, bounds the integrals on the remaining tails, giving
\eqref{eq:convex-uniform-BV-energy}.

The inner exponential in \eqref{eq:light-cone-superexponential} comes from the hyperbolic factors in the profile equations. On a nonconstant half-profile with limiting density $\mu_\infty$, let $b_\infty=p'(\mu_\infty)$. As $r=|\theta-\zeta|\to+\infty$, the coefficient of $q$ satisfies
\be
 \frac{q_r}{q}
 =-\frac{h(\mu_\infty)(1-b_\infty)}{2\eps}\,e^r(1+o(1)).
\ee
Indeed, $q\to0$ and $(p'(\mu)\cosh^2r-\sinh^2r)/\sinh r$ is asymptotic to $-(1-b_\infty)e^r/2$. Integrating the logarithmic derivative gives
\be
 \log|q(r)|=-\frac{h(\mu_\infty)(1-b_\infty)}{2\eps}\,e^r+o(e^r).
\ee
Exponentiating gives the decay of $q$ in $r$. Since $r=|\zeta|+\cO(1)$, the profile equations and integration from the endpoint yield, for some $a_\eps>0$,
\be
 |\mu-\mu_{l,r}|+|\mu'|+|\theta-\theta_{l,r}|+|\theta'|+|\theta''|
 \le C_\eps\exp\bigl(-a_\eps e^{|\zeta|}\bigr)
\ee
on the respective tails. Here one may decrease $a_\eps$ to absorb factors growing exponentially in $|\zeta|$ from the equations for the derivatives. The bound is faster than $e^{-N|\zeta|}$ for every fixed $N>0$. The form of this bound depends on the coordinate used to measure the approach to a light ray. In terms of $y=x/t$, the exact identity $e^{|\zeta|}=\sqrt{(1+|y|)/(1-|y|)}$ gives the bound
\be
 |\mu^\eps-\mu_{l,r}|+|\theta^\eps-\theta_{l,r}|
 \le C_\eps\exp\left(-\frac{a_\eps}{\sqrt{1-|y|}}\right)
\ee
near either ray, after adjusting the constants. In particular, the state differences vanish faster than every power of $1-|y|$. The asymptotic formula for $\log|q|$ in $r$ is invariant under boosts fixing the Riemann origin, since $\theta$ and $\zeta$ shift by the same constant, while $q$, $r$, and the limiting density are unchanged. Because the change of variable $y=\tanh\zeta$ introduces only factors of exponential growth when differentiating, the constant extension is $C^1$ in density and $C^2$ in rapidity at $y=\pm1$.
Inside the cone, the viscous flux is $\tau\del^\alpha U^\beta=S^\alpha(U^\beta)'$. The frame grows at most exponentially in $|\zeta|$, so this flux tends to zero at either light ray. The interior stress tends to the stress of the corresponding prescribed state, which is also the exterior stress. Integration by parts in the three regions of $t>0$ therefore produces no interface term: the normal traces of $T^{\alpha\beta}-\eps\tau_+\del^\alpha U^\beta$ agree across each light ray. The constant extension solves the conservative equation in the sense of distributions. Its bounded
perturbation from the Riemann initial data is supported in $|x|\le t$,
so the local $L^1$ error in the initial trace is $\cO(t)$.\qedhere
\end{enumerate}
\end{proof}

\subsection{Compactness and the initial trace}
\label{subsection:global-conditional-limit}

\paragraph{Compactness and endpoint states.}
The state and variation bounds give compactness on bounded intervals of rapidity. The exterior estimates recover the prescribed endpoint states $z_l$ and $z_r$ as $\zeta\to-\infty$ and $\zeta\to+\infty$, respectively, and give convergence on the whole rapidity line. These endpoint states determine the initial trace after the change back to spacetime coordinates.

\begin{lemma}[Compactness and recovery of the endpoint states]
\label{thm:global-conditional-limit}
Let $\eps_j\downarrow0$ and let $(\mu^j,\theta^j)$ be global profiles with fixed prescribed states $z_l,z_r$. Suppose their values stay in a fixed compact  subset where $h>0$ and $0<p'<1$, and every junction $\theta^j=\zeta$ satisfies $(\theta^j)'=0$. Then there is a limiting profile $(\mu,\theta)$ such that, along a subsequence, the differences $(\mu^j-\mu,\theta^j-\theta)$ tend to zero in $L^1(\mathbb R,d\zeta)$. This profile gives a weak solution of the Euler equations with bounded variation and the prescribed Riemann initial trace. It equals the corresponding endpoint states outside a bounded interval of rapidities, so its nonconstant region lies strictly inside the light cone.
\end{lemma}

\begin{proof}
\textbf{Local compactness and the Euler equations.} Boundedness of the rapidity gives opposite signs of $\theta^j-\zeta$ at the two ends. At every zero its derivative is $-1$, so there is one junction. Proposition~\ref{prop:nonlinear-monotonicity} bounds the variations by twice the ranges of the two state variables. Helly's theorem gives local strong convergence. For a compactly supported smooth vector test-function $\varphi(\zeta)$, integration of the vector profile equation gives
\be
-\int_{\mathbb R}T^{\alpha\beta}(\mu^j,U^j)(S_\alpha\varphi_\beta)'\,d\zeta=-\eps_j\int_{\mathbb R}(U^j)^ {\beta\prime}\varphi_\beta'\,d\zeta.
\ee
The right-hand side is bounded in absolute value by $\eps_j\|\varphi'\|_\infty\TV(U^j)$ and tends to zero. Strong local convergence on the compact state set gives strong local convergence of every stress component, yielding the Euler profile equation in the sense of distributions. To recover the spacetime equation, let $\chi_\beta(\tau,\zeta)$ be a smooth test function compactly supported in the interior cone and written in inertial components. The volume element is $\tau\,d\tau\,d\zeta$, while the coordinate gradients are $\del_\alpha\tau=-R_\alpha$ and $\del_\alpha\zeta=S_\alpha/\tau$. In the integral of $T^{\alpha\beta}\del_\alpha\chi_\beta$, integration by parts in $\tau$ therefore replaces the term $-\tau R_\alpha\del_\tau\chi_\beta$ by $R_\alpha\chi_\beta$. Since $S_\alpha'=R_\alpha$, the resulting integral at each $\tau$ is the integral of $T^{\alpha\beta}(S_\alpha\chi_\beta)'$, which vanishes by the Euler profile equation.

\textbf{Control of the exterior tails.}
To estimate the approach to the endpoint states uniformly, put $q^j=(\theta^j)'$ and choose $Z_0$ larger than the rapidity bound plus $\operatorname{arctanh}c^*+2$, with $c\le c^*<1$ on the state set. We write $\delta^j\coloneqq\theta^j-\zeta$ and $d_0\coloneqq\operatorname{arctanh}c^*+1$. On $\zeta\ge Z_0-1$ we have $\delta^j\le-d_0$, and on $\zeta\le-Z_0+1$ we have $\delta^j\ge d_0$. Rewriting \eqref{eq:A-coefficient-main} as
\be
A=h\frac{\cosh^2\delta^j}{\sinh\delta^j}
\bigl(\tanh^2\delta^j-c^2\bigr)
\ee
shows that $A$ has the sign of $\delta^j$ and obeys $|A|\ge a_*$, where $a_*\coloneqq\min h\,(\tanh^2 d_0-(c^*)^2)>0$ and the minimum is taken on the compact state set. The coefficient $B$ is bounded on these exterior intervals by $(c^*)^2\coth d_0$. Since $G'/G=-Bq^j$ and $\int|q^j|$ is uniformly bounded, ratios of the integrating factor are uniformly bounded on either exterior interval. Some $\eta_j\in[Z_0-1,Z_0-1/2]$ has $|q^j(\eta_j)|\le2\TV(\theta^j)$. The exact formula for the rapidity derivative therefore gives
\bse
\be
\label{eq:uniform-exterior-q}
 |q^j(\zeta)|\le C
 e^{-a_*(\zeta-Z_0+1/2)/\eps_j},\qquad \zeta\ge Z_0,
\ee
\be
 |\theta^j(\zeta)-\theta_r|+|\mu^j(\zeta)-\mu_r|
 \le C\eps_j e^{-a_*(\zeta-Z_0+1/2)/\eps_j}.
\ee
\ese
The second estimate follows by integration and the reconstruction bound $|\mu^{j\prime}|\le C|q^j|+C\eps_j|q^j|^2$. Backward integration from the analogous unit strip gives the left estimates. These arguments also apply to constant exterior branches.

The limiting profile is thus constant for $|\zeta|>Z_0$. Its constant extension across $x=\pm t$ satisfies the Euler equations without an interface jump. Its perturbation from the Riemann step is bounded and supported in $|x|\le t\tanh Z_0$, giving an $\cO(t)$ local $L^1$ error in the initial trace. The exterior estimates also upgrade the local convergence in $\zeta$ to strong convergence on the whole rapidity line: integration gives a bound $C\eps_j^2e^{-a_*/(2\eps_j)}$ for the two exterior $L^1(d\zeta)$ errors beyond $\pm Z_0$, while the limit equals the prescribed states there. On the fixed interval $[-Z_0,Z_0]$, boundedness and Helly convergence give strong $L^1$ convergence. Since $d(x/t)=\operatorname{sech}^2\zeta\,d\zeta$, strong convergence in the ray variable follows as well.
\end{proof}

\begin{proof}[Compactness conclusion for Theorem~\ref{thm:convex-global-profiles}]
The local junction construction and the global state estimates established above satisfy the assumptions of the preceding lemma. It gives a subsequential Euler limit and its Riemann trace, with $\mu\ge m_p$. Section~\ref{subsection:unique-euler-limit} below proves its entropy inequalities and identifies every such limit with the same classical Riemann solution, completing the proof of the theorem.
\end{proof}

\paragraph{Uniformity on compact sets of data.}
The construction permits the prescribed states to vary in a compact set with a viscosity threshold uniform on the compact set. The bounds below are also uniform over pressures with the same sound-speed bounds.

\begin{corollary}[Uniform estimates for compact sets of Riemann data]
\label{cor:compact-riemann-data}
Fix sound-speed bounds $0<b_-\le b_+<1$ and a nonempty compact set
\be
 K\Subset\big((0,+\infty)\times\mathbb R\big)^2
\ee
of pairs of Riemann states. There are constants $\eps_K,m_K,M_K,\Theta_K,C_K>0$ such that the existence and uniform estimates of Theorem~\ref{thm:convex-global-profiles} hold for every pair in $K$ and every $0<\eps<\eps_K$, for every pressure satisfying \eqref{eq:convex-pressure-assumptions} with these bounds. The constants depend only on $K,b_-,b_+$.
If the prescribed pairs converge in $K$, the pressures converge in $C^2$ on compact positive density intervals to a pressure satisfying \eqref{eq:convex-pressure-assumptions}, and $\eps_j\downarrow0$, every sequence of the constructed profiles converges strongly in $L^1_{\loc}(\mathbb R)$ to a bounded-variation Euler profile with the limiting Riemann data. Write $z^j$ for the viscous profiles and $z$ for the limiting Euler profile. If $S_j$ and $S$ are the steps at $\zeta=0$ with the prescribed and limiting endpoint states, respectively, then
\be
\label{eq:compact-data-centered-convergence}
(z^j-S_j)-(z-S)\longrightarrow 0
\qquad\text{in }L^1(\mathbb R,d\zeta).
\ee
The nonconstant part of the limiting profile lies in $|x/t|\le\tanh Z_K<1$, where $Z_K$ depends only on $K,b_-,b_+$.
\end{corollary}

\begin{proof}
We choose $0<a\le A$ and $T\ge1$ such that all endpoint densities lie in $[a,A]$ and all endpoint rapidities lie in $[-T,T]$. The paths \eqref{eq:convex-data-path} stay in the same box, so $|D|\le2T$ throughout each continuation. The pressure bounds imply $(1+b_-)\mu\le h(\mu)\le(1+b_+)\mu$. Choose $\bar\eps=b_-(1+b_-)a$ and define
\be
\begin{gathered}
 m_K\coloneqq a\exp(-L_{0,p}-2TL_{1,p}),\qquad
 \Theta_K\coloneqq T+\frac12\log\frac{(1+b_+)A}{(1+b_-)a},\\
 M_K\coloneqq\max\left\{A,
 \frac{(1+b_+)Ae^{4T}+4\bar\eps T e^{2T}}{2b_-}\right\}.
\end{gathered}
\ee
The constants $L_{0,p},L_{1,p}$ in \eqref{eq:convex-density-floor-constants} depend only on the pressure bounds. The formulas \eqref{eq:convex-upper-state-constants} and \eqref{eq:convex-uniform-floor} therefore give $m_K\le\mu\le M_K$ and $|\theta|\le\Theta_K$ for every matched preimage along every such path. We may take
\be
 \eps_K\coloneqq\frac12\min\{\bar\eps,b_-(1+b_-)m_K,(1+b_-)m_K/2\}>0.
\ee
This choice gives the threshold required in \eqref{eq:convex-epsilon-threshold} for each pressure and each data pair. The continuation argument can be applied separately to each pressure since its derivatives are bounded on the compact density interval. This threshold uses only the sound-speed bounds. Genuine nonlinearity is needed for each pressure to bound the positive rapidity derivative, but the value of that bound does not enter the viscosity threshold or the state estimates.

The variation bounds are $\TV(\mu)\le2(M_K-m_K)$ and $\TV(\theta)\le4\Theta_K$. The current identity used in the fourth part of the construction gives a uniform energy bound on a fixed interior interval, after choosing the normalization $n(m_K)=1$. The inequalities $1\le n(\mu)\le\exp((M_K-m_K)/((1+b_-)m_K))$ and $n/h\ge 1/((1+b_+)M_K)$ give common coefficient bounds. The exterior argument of Lemma~\ref{thm:global-conditional-limit} uses only these state and variation bounds and the strict upper bound on the sound speed. It gives uniform constants in \eqref{eq:uniform-exterior-q}, completing the energy estimate and fixing $Z_K$. Helly compactness gives strong local convergence of the states, and uniform convergence of the pressures on $[m_K,M_K]$ gives strong local convergence of the stresses. The profile equation therefore gives the Euler limit, while the same exterior estimates, measured relative to the varying endpoint states, identify its endpoint states with the limiting data. On $[-Z_K,Z_K]$, the convergence of the prescribed states gives $S_j\to S$ in $L^1$, so the local profile convergence also gives convergence after subtraction of these steps. Outside this interval $z-S=0$, and the exterior bounds on $z^j-S_j$ give a total $L^1(d\zeta)$ error at most $C\eps_j^2 e^{-a_*/(2\eps_j)}$, with common constants. The $C^2$ convergence of the pressures and strict positivity of $\mathcal G$ for the limiting pressure give a uniform positive lower bound for $\mathcal G_j$ on $[m_K,M_K]$. Thus the upper bounds for $q_j$ are uniform, and the entropy argument and classification in Section~\ref{subsection:unique-euler-limit} below apply to every cluster point with the limiting pressure and data. Thus the limit is independent of extraction, which proves \eqref{eq:compact-data-centered-convergence} along the whole sequence. The negative derivative bound and the amplitude bounds remain estimates at fixed $\eps$. The positive derivative bound is uniform in $\eps$ for a fixed pressure, or for the $C^2$-convergent pressure sequences just considered.
\end{proof}

\subsection{Entropy inequalities and identification of the Euler limit}
\label{subsection:unique-euler-limit}

\paragraph{Entropy production of the viscous profiles.}

The particle density $n$ defined above gives a direct entropy balance for the viscous profiles. In the interior of the light cone, the divergence of $nU$ is
\be
 \del_\alpha(nU^\alpha)
 =\frac{1}{\tau}\left[
 \bigl(n\sinh(\theta-\zeta)\bigr)'
 +n\cosh(\theta-\zeta)\right]
 =\frac{\eps n}{\tau h}(\theta')^2\ge0,
\ee
where the last equality follows from \eqref{eq:profile-current-identity}. Since the self-similar profile has $\del_\alpha\theta=(\theta'/\tau)S_\alpha$, the covector $\del\theta$ is spacelike wherever $\theta'\ne0$, and the production can also be written as $\eps\tau(n/h)\del_\alpha\theta\,\del^\alpha\theta$.

Outside the cone the current is constant on each side. Its two traces agree across each light ray by the endpoint estimates, so the distributional divergence has no interface term. Strong local convergence of the profiles on a compact  state subset gives strong local convergence of the currents on $t>0$. It follows that the limiting current has nonnegative divergence, or equivalently that the entropy pair given by the components of $-nU$ satisfies the entropy inequality. This argument uses the viscous production before taking the limit.

\paragraph{Convex entropy inequalities at shocks.}

For a general convex entropy pair we work with the limiting Euler profile in conservative variables. Its entropy balance has a diffuse part, which vanishes by the Euler equation, and a contribution at the jumps, whose sign follows from the variation of shock speed along the Hugoniot curve.

\begin{lemma}[Shock speeds under genuine nonlinearity]
\label{lemma:genuine-shock-speeds}
Fix a density $m>0$. For $\rho\ne m$ put $\Delta=\rho-m$, $h_l=m+p(m)$, and $a=(p(\rho)-p(m))/\Delta$. The squared shock speeds in the two fluid rest frames are
\be
S_l(\rho)=a\frac{h_l+\Delta}{h_l+a\Delta},\qquad
S_r(\rho)=a\frac{h_l+a\Delta}{h_l+\Delta}.
\ee
Their continuous extensions at $\rho=m$ equal $p'(m)$. The function $S_l$ is strictly increasing. For $\rho>m$ one has $S_l>p'(m)$ and $S_r<p'(\rho)$, with both inequalities reversed for $\rho<m$.
\end{lemma}
\begin{proof}
All denominators are positive, since $h_l+\Delta=\rho+p(m)$ and $h_l+a\Delta=m+p(\rho)$. Define $T(\rho)\coloneqq p'(\rho)-S_r(\rho)$. Taylor expansion at $m$ gives
\bse
\be
T(\rho)=\frac{\mathcal G(m)}{2h_l}(\rho-m)+o(|\rho-m|).
\ee
At a zero of $T$ with $\rho\ne m$, writing $b=p'(\rho)$ and using $a'=(b-a)/\Delta$ gives
\be
a'=\frac{a(a-1)}{h_l+\Delta},\qquad
S_r'=-\frac{2b(1-b)}{h(\rho)},\qquad
T'=\frac{\mathcal G(\rho)}{h(\rho)}>0.
\ee
Every zero away from $m$ would therefore be an upward crossing. The local signs furnished by the expansion exclude a first zero to the right of $m$ and a last zero to its left. Thus $T>0$ for $\rho>m$ and $T<0$ for $\rho<m$. Finally, direct differentiation gives
\be
S_l'(\rho)=\frac{h_l(h_l+\Delta)}{\Delta(h_l+a\Delta)^2}\,T(\rho)>0.
\ee
\ese
The strict inequalities follow from $S_l(m)=p'(m)$ and the sign of $T$.
\end{proof}

\begin{proposition}[Convex entropy inequalities for the Euler limit]
\label{prop:convex-entropy-limit}
Let $W=(T^{00},T^{01})$ be a limit in Theorem~\ref{thm:convex-global-profiles} or Corollary~\ref{cor:compact-riemann-data}, and write $F(W)=(T^{01},T^{11})$. Every $C^2$ entropy pair $(\mathcal E,\mathcal Q)$ on the conservative cone $\mathcal C$ with $\mathcal E$ convex and $D\mathcal Q=D\mathcal E\,DF$ satisfies
\be
\del_t\mathcal E(W)+\del_x\mathcal Q(W)\le0
\ee
in the sense of distributions on $t>0$.
\end{proposition}

\begin{proof}
  \textbf{1. \sl Orientation of the jumps.} The upper derivative bound in Lemma~\ref{lemma:convex-rapidity-derivative-bound} is independent of $\eps$ on the compact state set. Hence $D\theta^\eps\le Q_+\dd\zeta$ in the sense of distributions, and local strong convergence gives $D\theta\le Q_+\dd\zeta$. The atomic part is nonpositive, so at every jump $\theta_r\le\theta_l$, where the subscripts now denote the traces. If the two rapidities agree, a boost to their common rest frame gives the Rankine--Hugoniot relation $p(\mu_r)=p(\mu_l)$ from the momentum equation. The strict increase of $p$ then gives equality of the densities. Every nontrivial jump thus satisfies $\theta_r<\theta_l$.

 \textbf{2. \sl Entropy across a shock.}
We compute the Hugoniot curve from the left state in its rest frame. We write $m\coloneqq\mu_l$, $h_l\coloneqq m+p(m)$, and let $\rho>0$ denote the right density. We set $\Delta\coloneqq\rho-m$. If $v$ is the right velocity in this frame and $\sigma$ the shock speed, the two Rankine--Hugoniot equations read
\be
h(\rho)\frac{v}{1-v^2}
=\sigma\left(\Delta+h(\rho)\frac{v^2}{1-v^2}\right),
\qquad
p(\rho)-p(m)+h(\rho)\frac{v^2}{1-v^2}
=\sigma h(\rho)\frac{v}{1-v^2}.
\ee
Eliminating $\sigma$ and then using the first equation gives
\be
v^2=\frac{\Delta(p(\rho)-p(m))}{(m+p(\rho))(\rho+p(m))},
\qquad
E_r-m=\Delta\frac{m+p(\rho)}{h_l}.
\ee
The first identity excludes $\rho=m$ at a nontrivial jump. With $a\coloneqq(p(\rho)-p(m))/\Delta$, the condition $v<0$ selects the branch with
\be
\sigma=-\operatorname{sgn}(\Delta)
\sqrt{\frac{a(h_l+\Delta)}{h_l+a\Delta}}.
\ee
Conversely these formulas with $v<0$ give a Hugoniot arc from $\rho=m$ to $\rho=\mu_r$, taking its initial value by continuity. Its speed has magnitude less than one because $0<a<1$. Lemma~\ref{lemma:genuine-shock-speeds} shows that $S_l=\sigma^2$ increases strictly with $\rho$. On the arc with $\Delta>0$, the density increases and $\sigma=-\sqrt{S_l}$ decreases. On the arc with $\Delta<0$, the density decreases and $\sigma=\sqrt{S_l}$ also decreases. Boosting back preserves this order, since relativistic velocity addition is increasing on $(-1,1)$.

In the original coordinates, parametrize the arc by $\lambda\in[0,1]$ with $\rho=m+\lambda(\mu_r-m)$, and denote its state and speed by $W_\lambda$ and $s_\lambda$. Differentiating the Rankine--Hugoniot identity along the arc gives
\be
\frac{d}{d\lambda}\left(\mathcal Q(W_\lambda)-\mathcal Q(W_l)
-s_\lambda(\mathcal E(W_\lambda)-\mathcal E(W_l))\right)
=s_\lambda'\left(D\mathcal E(W_\lambda)(W_\lambda-W_l)
-\mathcal E(W_\lambda)+\mathcal E(W_l)\right).
\ee
The expression in parentheses on the right is nonnegative by convexity on $\mathcal C$, and $s_\lambda'\le0$. The entropy jump expression is zero at $\lambda=0$ and hence nonpositive at $\lambda=1$.

 \textbf{3. \sl The diffuse part of the entropy balance.}
Finally, in the ray variable $y=x/t$, the diffuse part of the profile derivative satisfies $(DF(W)-yI)D^{\rm d}W=0$ by the weak Euler equation and the BV chain rule~\cite{Volpert1967}. The diffuse part of $D\mathcal Q(W)-yD\mathcal E(W)$ consequently vanishes by entropy compatibility. At a jump at $y=s$, the atom of this measure is $[\mathcal Q]-s[\mathcal E]\le0$ by the preceding argument. The constant exterior regions introduce no additional atoms. To recover the spacetime inequality, take $\phi\in C_c^\infty((0,+\infty)\times\mathbb R)$ and set $\psi(t,y)=\phi(t,ty)$. The substitution $x=ty$ gives
\be
\begin{aligned}
\big\langle\del_t\mathcal E+\del_x\mathcal Q,\phi\big\rangle
 &=-\int_0^{+\infty}\!\int_{\mathbb R}
 \bigl(t\mathcal E(W)\del_t\psi+(\mathcal Q(W)-y\mathcal E(W))\del_y\psi\bigr)\,dy\,dt\\
 &=\int_0^{+\infty}\!\int_{\mathbb R}
 \psi(t,y)\,d\bigl(D\mathcal Q(W)-yD\mathcal E(W)\bigr)(y)\,dt.
\end{aligned}
\ee
Integration by parts in $t$ cancels the term $-\mathcal E(W)\,dy$ from the derivative of $-y\mathcal E(W)$. The resulting measure is nonpositive, which completes the proof of Proposition~\ref{prop:convex-entropy-limit}.\qedhere 
\end{proof}

\paragraph{Identification of the limits with the classical Riemann solution.}

We identify every limit with bounded variation of the constructed viscous profiles with the classical solution of the Euler Riemann problem, consisting of at most one wave of each characteristic family. The proof first derives this wave structure from the properties inherited by the BV limits and then uses the two wave curves to identify the common intermediate state.

\begin{enumerate}
\item \textbf{The junction and the characteristic families.} The pressure assumptions~\eqref{eq:convex-pressure-assumptions} include genuine nonlinearity. We write $v=\tanh\theta$ and denote the two characteristic speeds by
\be
\lambda_1=\frac{v-c}{1-vc},\qquad \lambda_2=\frac{v+c}{1+vc}.
\ee
We take a convergent subsequence and, after further extraction, let $\zeta_*$ be the limit of its junction coordinates. To the left of $\zeta_*$ the limiting profile satisfies $\theta\ge\zeta$, and to its right it satisfies $\theta\le\zeta$. Both inequalities follow by passage to the limit on the corresponding open intervals. The rapidity is monotone on each interval, since we can extract a subsequence with a fixed sign for the rapidity derivative on each branch. At a possible jump at the junction, the two relative rapidities have opposite weak signs. In the rest frame of this jump, continuity of the energy flux $h\sinh(\theta-\zeta_*)\cosh(\theta-\zeta_*)$ forces both relative rapidities to vanish. Continuity of momentum flux then gives equality of the pressures and states.

On the left, only the first characteristic family can carry diffuse variation, since $\lambda_2>y$ when $\theta\ge\zeta$. On the right, only the second family can carry it, since $\lambda_1<y$ when $\theta\le\zeta$.

\item \textbf{The Lax inequalities at a jump.} We consider the left side. A nontrivial jump there has positive relative fluid velocities in its rest frame, since they have nonnegative signs and vanishing energy flux would give identical states. In the rest frame of the left state, such a jump has negative speed and therefore $\rho>m$ in the Hugoniot formulas above. The squared shock speeds measured in the two fluid rest frames, denoted by $S_l$ and $S_r$, follow from the Hugoniot speed formula above by taking each trace in turn as the reference state. With $a$ the pressure difference quotient between the two densities, we obtain
\be
 S_l=a\frac{\rho+p(m)}{m+p(\rho)},\qquad
 S_r=a\frac{m+p(\rho)}{\rho+p(m)}.
\ee
Lemma~\ref{lemma:genuine-shock-speeds} gives $S_l>p'(m)$ and $S_r<p'(\rho)$ under genuine nonlinearity alone. In the shock rest frame the fluid velocities are positive, so the first characteristic speed is positive on the left and negative on the right. Transforming back gives $\lambda_1(W_l)>y>\lambda_1(W_r)$. Reflection gives the second-family Lax inequalities on the right.

\item \textbf{At most one compression wave on each branch.} On a branch of nonincreasing rapidity, the diffuse derivative vanishes. We use $D^{\rm d}$ for the diffuse part of the derivative in the ray variable. On the support of this derivative the Euler profile equation implies $\tanh(\theta-\zeta)=c$ on the left and $\tanh(\theta-\zeta)=-c$ on the right, since there $\lambda_i=y$. The projection of the Euler equation on the fluid velocity gives
\be
 \sinh(\theta-\zeta)\,D^{\rm d}\mu
 +h\cosh(\theta-\zeta)\,D^{\rm d}\theta=0 .
\ee
Thus $D^{\rm d}\theta=-(c/h)D^{\rm d}\mu$ on the left and $D^{\rm d}\theta=(c/h)D^{\rm d}\mu$ on the right, as equalities of measures. The corresponding characteristic integral curves satisfy
\be
d(\operatorname{arctanh}\lambda_i)
  =\left(1+\frac{h p''}{2p'(1-p')}\right)d\theta .
\ee
The factor is strictly positive. By the BV chain rule the diffuse derivative of $\lambda_i$ is therefore nonpositive when the rapidity is nonincreasing. Each shock also decreases $\lambda_i$ by the strict Lax inequalities just proved. Hence $\lambda_i-y$ is strictly decreasing in $y$, and its precise representative can vanish at at most one point. The diffuse state derivative is supported where $\lambda_i=y$ and has no atoms, so it vanishes. At a shock $\lambda_i-y$ goes from positive to negative, and strict decrease rules out another such shock on the same branch.

\item \textbf{A single rarefaction on an expanding branch.} On a branch of nondecreasing rapidity there are no jumps, since every nontrivial jump decreases the rapidity. The characteristic relations imply that $\theta+\int^\mu c(w)/h(w)\dd w$ is constant on the left and $\theta-\int^\mu c(w)/h(w)\dd w$ on the right, by the BV chain rule. Thus each continuous branch remains on a single characteristic curve. On that curve the characteristic speed is a strictly increasing function of rapidity. The support of the state variation is contained in the closed set where $\lambda_i=y$ within the branch. Any open gap between two points of this support would have a constant state and therefore the same characteristic speed at both edges, although the two edge values of $y$ are distinct. The support is thus an interval, and $\lambda_i=y$ there determines the rarefaction, including its edges, from the two endpoint states. A branch with no variation is constant.

\item \textbf{The intermediate state and convergence of the full family.} It remains to match the two waves. The left wave curve, written as a rapidity of the intermediate state as a function of its density, is strictly decreasing. Its rarefaction part has derivative $-c/h$. On its shock part the decrease in rapidity has magnitude $\operatorname{arctanh}\sqrt{V}$, with
\be
 V=\frac{\rho-m}{\rho+p(m)}
      \frac{p(\rho)-p(m)}{m+p(\rho)},\qquad \rho>m.
\ee
Both factors increase strictly in $\rho$. By reflection the wave curve from the right state is strictly increasing in the intermediate density. Their intersection is therefore unique whenever it exists. Figure~\ref{fig:riemann-wave-curves} illustrates the intersection for the data used in Section~\ref{section==7} below. Existence of an intersection follows from subsequential compactness and the preceding classification. The speeds and the rarefaction profiles are determined by the endpoint states of each wave, so every subsequence has the same limit. For fixed data, Lemma~\ref{thm:global-conditional-limit} gives subsequential compactness in $L^1(\mathbb R,d\zeta)$ for the differences of profiles. Since every cluster point agrees with this same Euler profile, the full family converges in the sense of \eqref{eq:convex-profile-convergence}.
\end{enumerate}

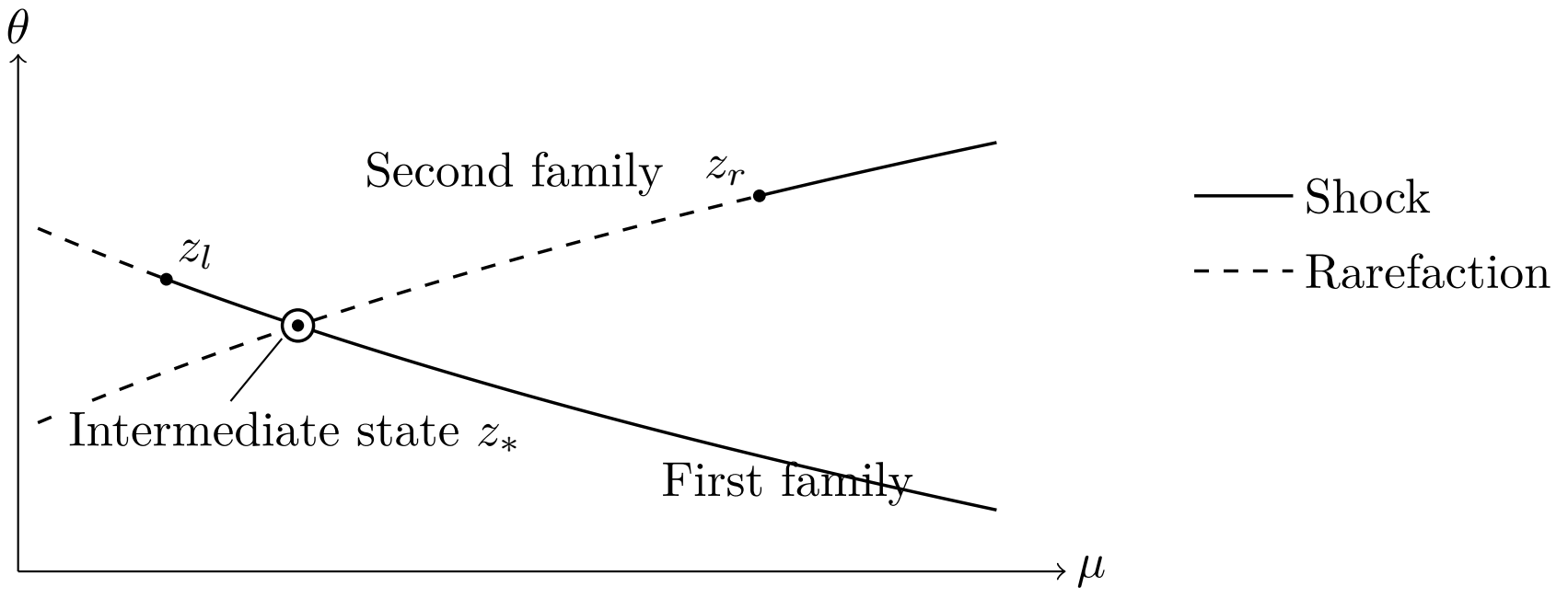
\begin{figure}[htbp]
\centering
\begin{tikzpicture}[x=7.1cm,y=6.0cm,
 every node/.style={font=\normalsize,inner sep=2pt},
 shock/.style={black,line width=.65pt}, rare/.style={shock,dashed}]
 \draw[->] (.85,-.35) -- (1.91,-.35) node[right] {$\mu$};
 \draw[->] (.85,-.35) -- (.85,.27) node[above] {$\theta$};
 \draw[rare] plot coordinates {(0.870000,0.061070) (0.872167,0.059985) (0.874333,0.058903) (0.876500,0.057823) (0.878667,0.056746) (0.880833,0.055671) (0.883000,0.054599) (0.885167,0.053529) (0.887333,0.052461) (0.889500,0.051396) (0.891667,0.050333) (0.893833,0.049273) (0.896000,0.048214) (0.898167,0.047159) (0.900333,0.046105) (0.902500,0.045054) (0.904667,0.044005) (0.906833,0.042959) (0.909000,0.041915) (0.911167,0.040873) (0.913333,0.039833) (0.915500,0.038795) (0.917667,0.037760) (0.919833,0.036727) (0.922000,0.035697) (0.924167,0.034668) (0.926333,0.033642) (0.928500,0.032618) (0.930667,0.031596) (0.932833,0.030576) (0.935000,0.029559) (0.937167,0.028544) (0.939333,0.027530) (0.941500,0.026519) (0.943667,0.025510) (0.945833,0.024504) (0.948000,0.023499) (0.950167,0.022496) (0.952333,0.021496) (0.954500,0.020498) (0.956667,0.019501) (0.958833,0.018507) (0.961000,0.017515) (0.963167,0.016525) (0.965333,0.015537) (0.967500,0.014551) (0.969667,0.013567) (0.971833,0.012585) (0.974000,0.011605) (0.976167,0.010627) (0.978333,0.009651) (0.980500,0.008677) (0.982667,0.007705) (0.984833,0.006735) (0.987000,0.005767) (0.989167,0.004801) (0.991333,0.003837) (0.993500,0.002875) (0.995667,0.001915) (0.997833,0.000956) (1.000000,0.000000)};
 \draw[shock] plot coordinates {(1.000000,0.000000) (1.014000,-0.006134) (1.028000,-0.012191) (1.042000,-0.018172) (1.056000,-0.024080) (1.070000,-0.029917) (1.084000,-0.035684) (1.098000,-0.041383) (1.112000,-0.047016) (1.126000,-0.052585) (1.140000,-0.058090) (1.154000,-0.063535) (1.168000,-0.068919) (1.182000,-0.074245) (1.196000,-0.079514) (1.210000,-0.084727) (1.224000,-0.089885) (1.238000,-0.094990) (1.252000,-0.100043) (1.266000,-0.105044) (1.280000,-0.109996) (1.294000,-0.114899) (1.308000,-0.119754) (1.322000,-0.124562) (1.336000,-0.129325) (1.350000,-0.134042) (1.364000,-0.138715) (1.378000,-0.143345) (1.392000,-0.147933) (1.406000,-0.152479) (1.420000,-0.156985) (1.434000,-0.161450) (1.448000,-0.165876) (1.462000,-0.170264) (1.476000,-0.174614) (1.490000,-0.178927) (1.504000,-0.183204) (1.518000,-0.187445) (1.532000,-0.191650) (1.546000,-0.195821) (1.560000,-0.199959) (1.574000,-0.204062) (1.588000,-0.208134) (1.602000,-0.212172) (1.616000,-0.216180) (1.630000,-0.220156) (1.644000,-0.224101) (1.658000,-0.228016) (1.672000,-0.231902) (1.686000,-0.235758) (1.700000,-0.239586) (1.714000,-0.243385) (1.728000,-0.247156) (1.742000,-0.250900) (1.756000,-0.254617) (1.770000,-0.258307) (1.784000,-0.261972) (1.798000,-0.265610) (1.812000,-0.269223) (1.826000,-0.272810) (1.840000,-0.276374)};
 \draw[rare] plot coordinates {(0.870000,-0.171837) (0.882167,-0.165778) (0.894333,-0.159796) (0.906500,-0.153887) (0.918667,-0.148051) (0.930833,-0.142285) (0.943000,-0.136588) (0.955167,-0.130958) (0.967333,-0.125394) (0.979500,-0.119894) (0.991667,-0.114456) (1.003833,-0.109080) (1.016000,-0.103763) (1.028167,-0.098505) (1.040333,-0.093304) (1.052500,-0.088158) (1.064667,-0.083068) (1.076833,-0.078030) (1.089000,-0.073046) (1.101167,-0.068112) (1.113333,-0.063229) (1.125500,-0.058395) (1.137667,-0.053609) (1.149833,-0.048870) (1.162000,-0.044177) (1.174167,-0.039530) (1.186333,-0.034928) (1.198500,-0.030369) (1.210667,-0.025853) (1.222833,-0.021379) (1.235000,-0.016946) (1.247167,-0.012554) (1.259333,-0.008202) (1.271500,-0.003888) (1.283667,0.000387) (1.295833,0.004625) (1.308000,0.008825) (1.320167,0.012989) (1.332333,0.017118) (1.344500,0.021211) (1.356667,0.025270) (1.368833,0.029295) (1.381000,0.033287) (1.393167,0.037246) (1.405333,0.041172) (1.417500,0.045067) (1.429667,0.048930) (1.441833,0.052763) (1.454000,0.056566) (1.466167,0.060339) (1.478333,0.064082) (1.490500,0.067797) (1.502667,0.071483) (1.514833,0.075141) (1.527000,0.078772) (1.539167,0.082375) (1.551333,0.085952) (1.563500,0.089503) (1.575667,0.093027) (1.587833,0.096526) (1.600000,0.100000)};
 \draw[shock] plot coordinates {(1.600000,0.100000) (1.604000,0.101137) (1.608000,0.102271) (1.612000,0.103402) (1.616000,0.104530) (1.620000,0.105656) (1.624000,0.106780) (1.628000,0.107901) (1.632000,0.109019) (1.636000,0.110134) (1.640000,0.111247) (1.644000,0.112358) (1.648000,0.113466) (1.652000,0.114571) (1.656000,0.115674) (1.660000,0.116774) (1.664000,0.117872) (1.668000,0.118968) (1.672000,0.120061) (1.676000,0.121151) (1.680000,0.122239) (1.684000,0.123325) (1.688000,0.124408) (1.692000,0.125488) (1.696000,0.126567) (1.700000,0.127643) (1.704000,0.128716) (1.708000,0.129787) (1.712000,0.130856) (1.716000,0.131923) (1.720000,0.132987) (1.724000,0.134049) (1.728000,0.135108) (1.732000,0.136165) (1.736000,0.137220) (1.740000,0.138273) (1.744000,0.139323) (1.748000,0.140371) (1.752000,0.141417) (1.756000,0.142461) (1.760000,0.143502) (1.764000,0.144541) (1.768000,0.145578) (1.772000,0.146613) (1.776000,0.147645) (1.780000,0.148676) (1.784000,0.149704) (1.788000,0.150730) (1.792000,0.151754) (1.796000,0.152775) (1.800000,0.153795) (1.804000,0.154812) (1.808000,0.155828) (1.812000,0.156841) (1.816000,0.157852) (1.820000,0.158861) (1.824000,0.159868) (1.828000,0.160873) (1.832000,0.161876) (1.836000,0.162877) (1.840000,0.163875)};
 \fill (1,0) circle (1.3pt) node[above right] {$z_l$};
 \fill (1.6,.1) circle (1.3pt) node[above left] {$z_r$};
 \node[anchor=west] at (1.49,-.245) {First family};
 \node[anchor=west] at (1.19,.126) {Second family};
 \filldraw[fill=white,draw=black,line width=.65pt] (1.1331187,-0.0553921) circle (3.2pt);
 \fill (1.1331187,-0.0553921) circle (1.25pt);
 \draw[thin] (1.117,-.071) -- (1.065,-.146);
 \node[anchor=north west,fill=white] at (.89,-.146) {Intermediate state $z_*$};
 \draw[shock] (2.04,.10) -- (2.14,.10) node[right] {Shock};
 \draw[rare] (2.04,.01) -- (2.14,.01) node[right] {Rarefaction};
\end{tikzpicture}
\caption{Wave curves in the density--rapidity plane for $z_l=(1,0)$ and $z_r=(1.6,0.1)$, with the pressure law~\eqref{eq:convex-nonaffine-example} below. The descending curve connects $z_l$ to a possible intermediate state by a first-family wave. The ascending curve connects a possible intermediate state to $z_r$ by a second-family wave. Their intersection $z_*$ selects a first-family shock and a second-family rarefaction.}
\label{fig:riemann-wave-curves}
\end{figure}
\FloatBarrier

\section{Travelling waves near an isolated shock}
\label{section:traveling-waves}

\subsection{Travelling-wave equations}

We derive a scalar traveling-wave equation by rescaling at a shock. We then justify the rescaling limit for the viscous profiles and use the limiting wave to determine the width of the compression layer.

\begin{enumerate}
\item \textbf{Reduction to a scalar equation.}
In this section, the subscripts $l,r$ denote the two traces at the selected shock, rather than necessarily the original Riemann states. We consider a discontinuity of an Euler limit at cone rapidity $\zeta_s$ and introduce the coordinate $\eta=(\zeta-\zeta_s)/\eps$. On this scale the vector profile equation becomes
\be
S_\alpha(\zeta_s+\eps\eta)\frac{dT^{\alpha\beta}}{d\eta}
=\frac{d^2U^\beta}{d\eta^2}.
\ee
Freezing the frame at $\zeta_s$ gives an autonomous traveling-wave equation, for which we introduce the vectors and relative rapidity
\be
R_s\coloneqq(\cosh\zeta_s,\sinh\zeta_s),\quad S_s\coloneqq(\sinh\zeta_s,\cosh\zeta_s),\quad \delta\coloneqq\theta-\zeta_s.
\ee
Integration of the autonomous equation from the left state gives
\bse
\be
\label{eq:traveling-vector}
U_\eta=S_{s,\alpha}T^{\alpha\cdot}(\mu,U)-F_l,
\qquad F_l\coloneqq S_{s,\alpha}T^{\alpha\cdot}(\mu_l,U_l)
=A R_s+B S_s.
\ee
Since $U\cdot U_\eta=0$, projection on $U$ and on $N$ gives
\be
\label{eq:traveling-scalar}
\mu=A\coth\delta-B,
\qquad \delta_\eta=A\sinh\delta+(p(\mu)-B)\cosh\delta,
\ee
\ese
on an interval where $\sinh\delta\ne0$. With $u\coloneqq\tanh\delta$, these equations reduce to
\be
\label{eq:traveling-general}
\mu=\frac Au-B,\qquad
u_\eta=\sqrt{1-u^2}\,H(u),\qquad
H(u)\coloneqq p(A/u-B)+Au-B.
\ee
The Rankine--Hugoniot relations impose the same flux at both endpoints. For states satisfying the density relation in \eqref{eq:traveling-general}, equality of their fluxes to $F_l$ is equivalent to the two conditions $H(u_l)=H(u_r)=0$.

\item \textbf{Monotonicity and exponential tails.}
For a first-family compression wave we have $A>0$ and $0<u_r<u_l<1$. At an endpoint the flux relation gives
\be
H'(u)=A\left(1-\frac{p'(\mu)}{u^2}\right).
\ee
The Lax inequalities from Lemma~\ref{lemma:genuine-shock-speeds} imply $H'(u_l)>0$ and $H'(u_r)<0$. An additional zero between $u_r$ and $u_l$ would determine another density strictly between $\mu_l$ and $\mu_r$ with the same fluxes and the same squared shock speed $S_l=u_l^2$ in the left fluid rest frame. Strict monotonicity of $S_l$ excludes such a zero. Therefore $H<0$ on $(u_r,u_l)$ and its two endpoint zeros are simple. Separation of variables gives a unique decreasing connection up to translation, with increasing density and exponential tails. The linearized rates are $\sqrt{1-u_l^2}H'(u_l)>0$ on the left and $\sqrt{1-u_r^2}H'(u_r)<0$ on the right. This argument uses genuine nonlinearity, without requiring convexity of $H$ or of the pressure.

\end{enumerate}

For the isothermal law $p(\mu)=k^2\mu$, the two roots of $H$ give the factorization
\be
H(u)=\frac Au(u-u_l)(u-u_r),\qquad u_lu_r=k^2.
\ee
Substitution in \eqref{eq:traveling-general} determines the traveling wave by separation of variables.

\subsection{Localization and centering of the viscous layer}
\label{subsection:centered-shock-layers}

Suppose a sequence of the constructed profiles converges locally in $L^1$ to an Euler profile with a nontrivial isolated shock at $\zeta_s$. The energy flux locates the shock on the appropriate viscous branch. We treat positive flux first and return to negative flux in Section~\ref{subsection:shock-scale-convergence} below.

\begin{enumerate}
\item \textbf{The first-family shock.} The two shock-frame velocities have the same sign, since the Rankine--Hugoniot relations give the common energy flux
\be
A=\frac{h(\mu_l)u_l}{1-u_l^2}=\frac{h(\mu_r)u_r}{1-u_r^2}.
\ee
If this flux vanishes, both relative velocities are zero and equality of the momentum fluxes gives $p(\mu_l)=p(\mu_r)$. Strict pressure monotonicity then excludes a nontrivial jump. Proposition~\ref{prop:convex-entropy-limit} and its proof give $\theta_r<\theta_l$, so $u_r<u_l$. Consider first the positive-flux case, in which $0<u_r<u_l<1$. At an endpoint the derivative in \eqref{eq:traveling-general} is $H'(u)=A(1-c^2/u^2)$, by the flux identity above. The signs of $H'$ at its two roots thus give $u_l>c(\mu_l)$ and $u_r<c(\mu_r)$, identifying a first-family shock.

\item \textbf{Separation from the junction.} We choose continuity points $a<\zeta_s<b$ sufficiently close to the shock that $\theta(a)>\theta(b)$ and $\theta(b)-b>0$. The classification in Section~\ref{subsection:unique-euler-limit} places this shock strictly to the left of the limiting junction. We choose $a,b$ on that same side, so for sufficiently small $\eps$ they lie on one monotone viscous branch. Monotonicity of $\theta^\eps$ and local $L^1$ convergence give convergence at these continuity points, by bracketing each value between averages on adjacent intervals. Hence $\theta^\eps(b)-b\ge d>0$ and $\theta^\eps(a)>\theta^\eps(b)$ for some fixed $d>0$. The derivative has a fixed sign on this branch, and this ordering makes it negative. Consequently $(\theta^\eps-\zeta)'=q-1<0$ on $[a,b]$, so the relative rapidity is bounded below there by $d$. Both state variables are consequently monotone on this fixed interval.

\item \textbf{Centering by the rapidity midpoint.} We center the layer at $\zeta_\eps$ by requiring
\be
\theta^\eps(\zeta_\eps)-\zeta_s
=\tfrac12\bigl(\theta_l+\theta_r\bigr)-\zeta_s.
\ee
We choose continuity points to the left and right of $\zeta_s$ whose limiting values bracket the midpoint strictly. For small $\eps$, the same strict inequalities hold for the profiles, which therefore cross the midpoint. Such bracketing points can be taken arbitrarily close to $\zeta_s$, giving $\zeta_\eps\to\zeta_s$. The compact uniform bounds and the variation bounds persist under the change $\eta=(\zeta-\zeta_\eps)/\eps$.

\end{enumerate}

\subsection{Convergence on the shock scale}
\label{subsection:shock-scale-convergence}

The integrated flux gives local uniform compactness of the rescaled fields and identifies the limiting equation. Monotonicity and the midpoint condition then select the unique traveling wave.

\begin{enumerate}
\item \textbf{Compactness of the integrated flux.} We set
\be
P_\eps^\beta(\zeta)\coloneqq
S_\alpha(\zeta)T^{\alpha\beta}(\mu^\eps,U^\eps)-\eps(U^\eps)^{\beta\prime}.
\ee
The identity $P_\eps'=R_\alpha T^{\alpha\cdot}$ gives a uniform local Lipschitz bound. A bound for $P_\eps$ itself follows by averaging on any fixed interval: the stress is bounded and $\eps\int|U'|$ is bounded. After extraction, $P_\eps$ converges locally uniformly. Since $\eps U'\to0$ in the sense of distributions, its limit equals $S_\alpha T^{\alpha\cdot}$ away from the shock; continuity and Rankine--Hugoniot give $P_0(\zeta_s)=F_l$.

\item \textbf{The autonomous equation.} For the centered fields $V^\eps(\eta)\coloneqq U^\eps(\zeta_\eps+\eps\eta)$, the exact identity is
\be
V^\eps_\eta=S_\alpha(\zeta_\eps+\eps\eta)T^{\alpha\cdot}-P_\eps(\zeta_\eps+\eps\eta).
\ee
It bounds $V^\eps_\eta$ locally uniformly, so the velocities are locally equicontinuous. At $\zeta=\zeta_\eps+\eps\eta$, we decompose the exact flux in the moving frame as $P_\eps=A_\eps R(\zeta)+B_\eps S(\zeta)$. Orthogonality to $V^\eps$ gives
\be
\mu^\eps(\zeta)=A_\eps(\zeta)\coth(\theta^\eps(\zeta)-\zeta)-B_\eps(\zeta).
\ee
The preceding separation estimate bounds $\sinh(\theta^\eps-\zeta)$ away from zero. Uniform convergence of $P_\eps$ and of the velocities therefore gives local uniform convergence of the density as well. The exact equation then gives local uniform convergence of $V^\eps_\eta$ and the autonomous equation~\eqref{eq:traveling-vector}.

\item \textbf{Identification of the connection.} Monotonicity brackets the centered profiles by their values at fixed points on either side of the shock. Passing first to the zoom limit and then letting those points tend to $\zeta_s$ places the limiting relative velocity between $u_r$ and $u_l$. The limiting relative velocity is nonincreasing, and its limits at $\eta=\pm\infty$ must be roots of $H$; the absence of any other zero between the endpoints, proved for the scalar connection above, leaves only these two roots. The centering condition excludes either constant solution. The limit is therefore the unique normalized connection. For fixed Riemann data the outer profile is unique, so this identification gives convergence along the full family near each of its isolated shocks.

\item \textbf{The second-family shock.} For a negative energy flux, reflect $x$ to $-x$, which sends $\zeta$ to $-\zeta$ and the state to $(\mu,-\theta)$. The equations are invariant under this reflection, which interchanges the left and right traces and sends $(u_l,u_r)$ to $(-u_r,-u_l)$. The reflected shock has positive flux and is of the first family, so the original shock is of the second family. Applying the preceding argument to the reflected profiles gives $\theta^\eps-\zeta\le-d$ and $(\theta^\eps)'<0$ on a fixed interval about the original shock. The midpoint centering is preserved under reflection, with the sign of the stretched coordinate reversed. Reflecting the limiting traveling wave back gives the claimed local convergence for the second family.

\end{enumerate}

\subsection{Width of the compression layer}
The convergence on the shock scale determines the distance between two fixed interior density levels. Converting those crossing points to the ray variable gives the leading width of the layer.

\begin{corollary}[Shock-layer width]
\label{cor:shock-layer-width}
Along a sequence $\eps_j\downarrow0$ of viscosity parameters for which the profiles converge to an Euler profile with an isolated nontrivial shock, as in Sections~\ref{subsection:centered-shock-layers} and~\ref{subsection:shock-scale-convergence}, let $s=\tanh\zeta_s$ and fix $0<a<b<1$. For $\eps=\eps_j$ sufficiently small, let $y_a^\eps,y_b^\eps$ be the two ray coordinates in the layer at which the density takes the values $\mu_l+a(\mu_r-\mu_l)$ and $\mu_l+b(\mu_r-\mu_l)$, respectively. Denote by $\eta_a,\eta_b$ the positions of these same levels on the centered traveling wave. Then
\be
\label{eq:shock-layer-width}
\lim_{j\to\infty}\frac{|y_b^{\eps_j}-y_a^{\eps_j}|}{\eps_j}
=(1-s^2)|\eta_b-\eta_a|.
\ee
The same limit holds if the two densities used to define the levels in the viscous profile converge to $\mu_l$ and $\mu_r$, respectively.
\end{corollary}

\begin{proof} The density of the autonomous wave is strictly monotone between its two limiting values. For each specified level, we choose two points on either side of its crossing in the stretched coordinate $\eta$. Local uniform convergence preserves the strict bracketing inequalities for sufficiently small $\eps$, including when the prescribed levels vary and converge. The viscous density is strictly monotone near the shock, so the crossing is unique and converges in the stretched coordinate to the crossing on the wave.
  In cone rapidity, the two crossings are therefore at $\zeta_\eps+\eps(\eta_a+o(1))$ and $\zeta_\eps+\eps(\eta_b+o(1))$. Since $\zeta_\eps\to\zeta_s$ and $y=\tanh\zeta$, the mean value theorem gives \eqref{eq:shock-layer-width}.\qedhere
\end{proof}

\section{Numerical experiments}
\label{section==7}

\subsection{A shock and a rarefaction}

\paragraph{Viscous profiles and Euler reference.}
We compare the viscous profiles with the Euler solution, before testing the shock width for large data, finite propagation, and covariance under a Lorentz boost.
For this example we take the convex pressure
\be
\label{eq:convex-nonaffine-example}
p(\mu)\coloneqq a\mu+b\left[\mu-\mu_0\log\left(1+\frac{\mu}{\mu_0}\right)\right],
\qquad (a,b,\mu_0)=(0.15,0.30,1).
\ee
Its derivative increases from $a$ to $a+b$, so the sound speed satisfies the bounds of the main theorem.

\begin{figure}[htbp]
\centering\includegraphics[width=\textwidth]{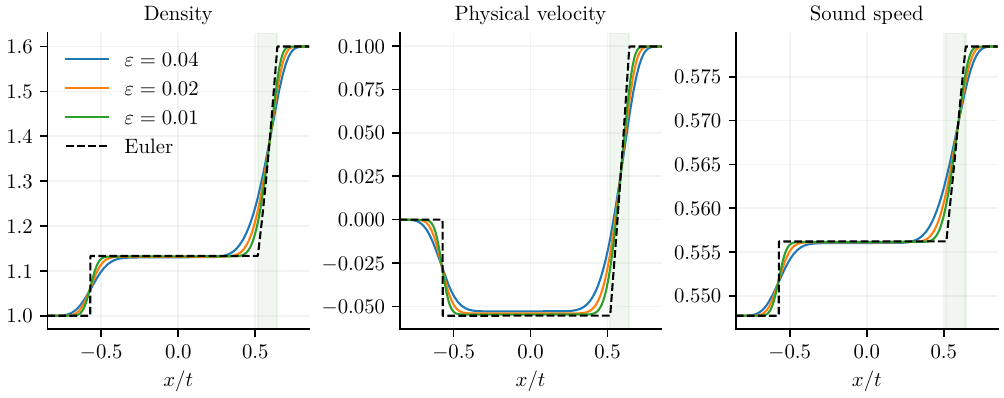}
\caption{Two-wave Riemann profiles for a nonaffine pressure: density, physical velocity, and variable sound speed. Dashed curves are the Euler reference; shading marks its second-family rarefaction fan.}
\label{fig:nonaffine-profiles}
\end{figure}

The states are $z_l=(1,0)$ and $z_r=(1.6,0.1)$ in density--rapidity coordinates. We integrate the half-profile equations \eqref{eq:fuchsian-fixed-epsilon} and solve the three matching equations for the junction density and the two amplitudes. The linear profile with coefficients fixed at the junction is factored out before integration. Adaptive fourth-order Runge--Kutta integration and damped Newton iteration give the profiles in Figure~\ref{fig:nonaffine-profiles}.

The Euler reference consists of a first-family Lax shock of speed $-0.57164$ and a second-family rarefaction with edge speeds $0.51677$ and $0.64116$. We compute the intermediate state from the Rankine--Hugoniot relations and the second-family Riemann invariant, with errors below $10^{-12}$. Figure~\ref{fig:nonaffine-profiles} compares the viscous profiles with this reference.

Refinement halves the starting distance from the junction, extends the integration interval from $r=4$ to $r=5$, and reduces the tolerance from $2\times10^{-7}$ to $2\times10^{-9}$. For $\eps=0.04,0.02,0.01$, the matching errors stay below $9\times10^{-10}$ and the differences in $(\mu,\theta)$ on $[-2,2]$ stay below $9.2\times10^{-4}$. The errors in the weighted flux identities are below $10^{-9}$.

\paragraph{Resolution of the compression layer.}

\begin{figure}[htbp]
\centering
\includegraphics[width=0.73\textwidth]{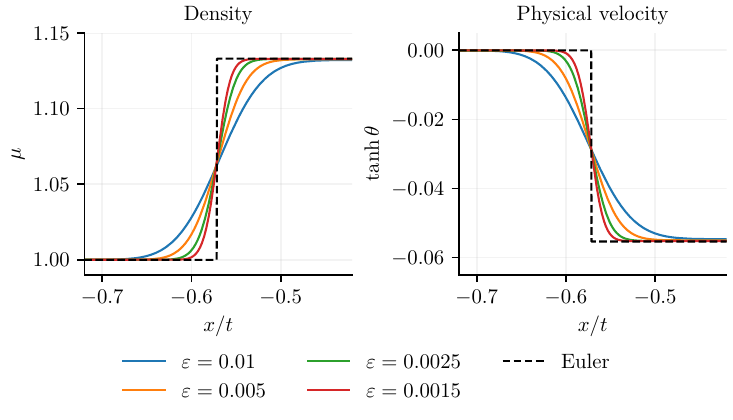}
\caption{Enlargement of the first-family compression wave in the same numerical profiles, on $-0.72\le x/t\le-0.42$. Density and physical velocity are shown for $\eps=0.01,0.005,0.0025,0.0015$ on the same fixed ray interval; dashed curves give the discontinuous Euler reference, whose shock speed is approximately $-0.57164$. The computed transitions become narrower as $\eps$ decreases. }
\label{fig:shock-zoom-numerical}
\end{figure}

Figure~\ref{fig:shock-zoom-numerical} resolves the compression layer for smaller viscosities on a fixed interval of $x/t$, so its narrowing can be compared directly with the Euler shock.

The compression in these data is weak, so its transition remains broad at $\eps=0.01$. For $\eps=0.0015$, the width between $10\%$ and $90\%$ of the density jump is $0.03087$ in $x/t$. The original and refined adaptive integrations place $37$ and $58$ points in this interval and give widths $0.030886$ and $0.030868$. Reducing the tolerance from $2\times10^{-10}$ to $2\times10^{-11}$ and the maximum step in $r$ from $0.015$ to $0.003$, moving the initial radius from $0.00125$ to $0.000625$, and extending the outer radius from $5.5$ to $6$ changes the plotted state by at most $2.7\times10^{-5}$ on the displayed interval. The refined endpoint error is below $10^{-7}$. These checks concern the resolved finite-$\eps$ layers; the autonomous traveling-wave figure uses a stretched coordinate and therefore displays a smooth transition.

\subsection{Compression for large Riemann data}

\begin{figure}[htbp]
\centering\includegraphics[width=0.73\textwidth]{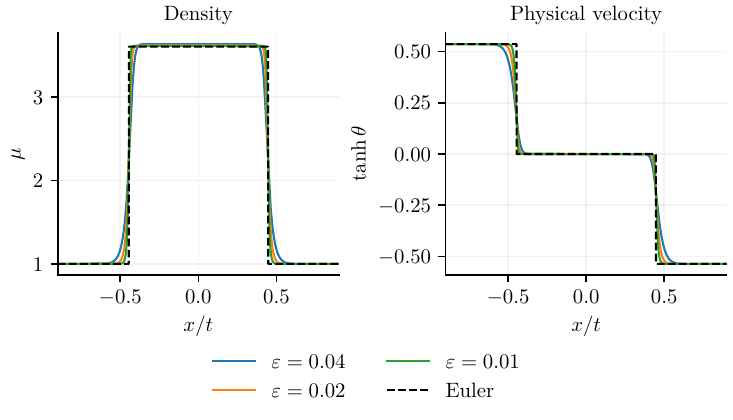}
\caption{Large-jump Riemann profiles for the convex nonaffine pressure, with endpoint states $(\mu_l,\theta_l)=(1,0.6)$ and $(\mu_r,\theta_r)=(1,-0.6)$. Density and velocity are plotted on the same ray interval for $\eps=0.04,0.02,0.01$. Dashed curves show the two-shock Euler reference with central state $(3.6041967360,0)$ and shock speeds $\pm0.4462285040$.}
\label{fig:large-jump-riemann}
\end{figure}

\paragraph{Symmetric Riemann data.}
For a large Riemann jump, we test convergence of the centered viscous layers to the traveling waves. We also compare the measured width divided by $\eps$ with the asymptotic coefficient in Corollary~\ref{cor:shock-layer-width}.
For the same nonaffine pressure, we next prescribe $z_l=(1,0.6)$ and $z_r=(1,-0.6)$. The velocities are approximately $0.53705$ and $-0.53705$, and the relative velocity has magnitude $\tanh(1.2)\simeq0.83365$. These colliding states produce two compression layers, with a central density more than three times the endpoint density. The symmetric Euler reference has zero intermediate rapidity, intermediate density $\mu_m\simeq3.6041967360$ and two shocks of speeds approximately $\pm0.4462285040$. We compute $\mu_m$ from
\be
 \tanh^2(0.6)=
 \frac{(p(\mu_m)-p(1))(\mu_m-1)}{(1+p(\mu_m))(\mu_m+p(1))},
\ee
and recover the shock speeds from the Rankine--Hugoniot equations, whose errors are below $10^{-14}$. Figure~\ref{fig:large-jump-riemann} compares this reference with the viscous profiles.

The refined integration uses a maximum step $0.001$ in $r$, tolerance $2\times10^{-12}$, initial radius $0.000625$, and outer radius $6$. The left compression layer, measured between $10\%$ and $90\%$ of its density rise, contains respectively $527$, $486$, and $462$ accepted integration points. Its width decreases from $0.07616$ to $0.04180$ and $0.02211$ as $\eps$ decreases through the three values. Comparison with the preceding integrations, whose maximum step is $0.015$, gives the differences in Table~\ref{table:large-jump-refinement}; these include interpolation error. Matching errors are below $4\times10^{-9}$ and the errors in the weighted flux identities are below $2\times10^{-13}$. The measured $L^1(d(x/t))$ errors in density and rapidity against the Euler reference decrease from $0.1836$ to $0.1011$ and $0.05354$. This finite computation illustrates large compression and resolved narrowing of both layers, without asserting a general convergence rate.

\begin{table}[htbp]\centering
\begin{tabular}{ccccc}\hline
$\eps$ & $\mu_*^\eps$ & $\TV(\mu)+\TV(\theta)$ & Layer width & Refinement difference\\\hline
0.04 & 3.63567 & 6.47133 & 0.07616 & $5.58\times10^{-4}$\\
0.02 & 3.61944 & 6.43888 & 0.04180 & $2.98\times10^{-4}$\\
0.01 & 3.61164 & 6.42328 & 0.02211 & $2.84\times10^{-4}$\\\hline
\end{tabular}
\caption{Large-jump compression: central density, total variation, width in $x/t$ between $10\%$ and $90\%$ of the left density rise, and maximum difference in $(\mu,\theta)$ under refinement on $[-0.999,0.999]$.}
\label{table:large-jump-refinement}
\end{table}

\paragraph{Continuation and comparison with traveling waves.}
We continue the matching problem from $\eps=0.01$ toward smaller viscosities, adapting the method of Schecter, Plohr, and Marchesin~\cite{SchecterPlohrMarchesin2004}. Symmetry reduces the calculation to one half-profile, with unknown junction density and logarithmic amplitude. A predictor based on the preceding two solutions extrapolates these quantities in $\eps$ and $1/\eps$, respectively, before Newton correction enforces the endpoint states.

We center each viscous layer and the corresponding traveling wave at rapidity $0.3$ for the left shock and $-0.3$ for the right shock. Density and velocity in the Euler shock frame are compared at spacing $0.001$ on $[-4,4]$ in $\eta=(\zeta-\zeta_\eps)/\eps$. Reflection reverses $\eta$ and the relative velocity, giving identical errors and widths for the two shocks.

\begin{figure}[htbp]
\centering
\includegraphics[width=0.66\textwidth]{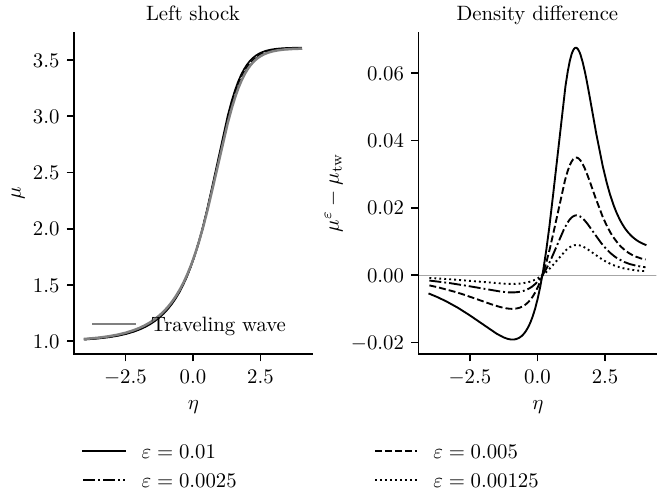}
\caption{Density and signed density difference for the left compression layer, centered at the midpoint of the limiting rapidities and plotted in the stretched coordinate $\eta$. The difference is the viscous density minus the autonomous traveling-wave density. Reflection gives the corresponding curves for the right shock with $\eta$ replaced by $-\eta$.}
\label{fig:left-shock-layer-comparison}
\end{figure}

The finer calculation uses a maximum step $\eps/20$ in $r$, tolerance $3\times10^{-13}$, and the interval $[0.0003125,4.5]$. We repeat the matching with maximum step $\eps/8$, tolerance $3\times10^{-11}$, and the interval $[0.000625,3.5]$. Centering and interpolating the two integrations separately gives differences below $8\times10^{-5}$ in either field on the stretched window. Halving the step in the autonomous velocity equation changes its values on the same grid by less than $10^{-10}$. The widths in Table~\ref{table:symmetric-continuation} are measured relative to the computed central density and the prescribed exterior density. Over the tested values, the computed central density approaches the Euler intermediate density. For the autonomous wave in \eqref{eq:traveling-general}, quadrature between the $10\%$ and $90\%$ density levels gives $(1-s^2)|\eta_{90}-\eta_{10}|\simeq2.3603$. The measured width divided by $\eps$ reaches $2.340$ at the smallest $\eps$, within $0.9\%$ of the limit predicted by Corollary~\ref{cor:shock-layer-width}. Figure~\ref{fig:left-shock-layer-comparison} shows where the density differs from the autonomous wave.

\begin{table}[htbp]
\centering
\begin{tabular}{rrrrrr}
\hline
$\eps$ & Width & Width$/\eps$ & Points & Density error & Velocity error \\
\hline
0.01 & 0.022106 & 2.211 & 925 & 0.06756 & 0.006804 \\
0.005 & 0.011405 & 2.281 & 1798 & 0.03492 & 0.003539 \\
0.0025 & 0.005798 & 2.319 & 3543 & 0.01778 & 0.001807 \\
0.00125 & 0.002924 & 2.340 & 7030 & 0.00897 & 0.000913 \\
\hline
\end{tabular}
\caption{Continuation at fixed symmetric Riemann data. The widths in $x/t$ are measured between $10\%$ and $90\%$ of the density change. The point counts refer to accepted integration points within each layer. The last two columns give the maximum sampled differences from the traveling waves in density and relative velocity. Reflection gives identical values for the two shocks.}
\label{table:symmetric-continuation}
\end{table}

\subsection{Finite propagation in spacetime}

\paragraph{First-order formulation.}
This experiment illustrates that an initially compact perturbation stays inside the region reached at light speed. Small numerical tails beyond this causal region arise from the discretization and diminish under mesh refinement.
To display propagation, we also evolve smooth compactly supported perturbations of a constant state using the conservative equation \eqref{eq:conservative-tau-model-main}. This time-dependent calculation includes the derivatives of $\tau$ inside the divergence. On a profile depending only on the cone rapidity these terms vanish, as observed in Section~\ref{subsection:self-similar-velocity-viscosity}.

We set $w\coloneqq\sinh\theta$, $z\coloneqq w_t$, $q\coloneqq w_x$, and $v\coloneqq w/\sqrt{1+w^2}$. Projection of the spacetime equation on $U$ gives $U^\alpha\del_\alpha\mu$ as its only derivative of the first-order variables, since $U\cdot\del_\alpha U=0$ and $U\cdot\Box U=-\del_\alpha U\cdot\del^\alpha U$. The spatial component of the equation contains the wave operator on $w$, together with first derivatives of the density. After solving for the time derivatives, the principal equations have the form
\be
\mu_t+v\mu_x=\text{lower order terms},\qquad
z_t-q_x=B\mu_x+\text{lower order terms},\qquad
q_t-z_x=0,
\ee
The coefficient is $B=-p'(\mu)/(\eps\tau)$. Indeed, the density derivatives in the spatial component of the stress divergence are $h'(\mu)\sqrt{1+w^2}w\mu_t+(h'(\mu)w^2+p'(\mu))\mu_x$. Substitution of the transport equation cancels the terms proportional to $h'(\mu)w^2$, and the coefficient of $z_t-q_x$ in the viscous term is $-\eps\tau$. Here the lower order terms contain no derivatives of $(\mu,z,q,w)$. The principal matrix is block triangular, with the scalar speed $v$ and the two wave speeds $-1$ and $1$. The equation $w_t=z$ adds a zero speed and no spatial derivatives. Since $|v|<1$, the density transport speed is distinct from the two wave speeds, and the principal matrix is smoothly diagonalizable on compact state subsets. The domain of dependence for smooth solutions is therefore bounded by the light speeds. For this experiment $\eps=0.03$ and the initial time is $t_0=2$. With
\bse
\be
\psi(s)=\begin{cases}\exp(1-(1-s^2)^{-1}),&|s|<1,\\0,&|s|\ge1,\end{cases}
\ee
we prescribe
\be
\mu(t_0,x)=1+0.35\psi((x+0.08)/0.15),\quad
 w(t_0,x)=0.12\psi((x-0.08)/0.15),\quad z(t_0,x)=0,
\ee
\ese
and impose the discrete compatibility $q=D_xw$. These perturbations are supported in $[-0.3,0.3]$. For a smooth solution, finite propagation therefore confines the support of $(\mu-1,z,q,w)$ to the compact interval $[-0.3-(t-t_0),0.3+(t-t_0)]$. At the final time $t=t_0+0.4$, the exact solution is the constant background outside $[-0.7,0.7]$. This support statement concerns the perturbation, rather than the positive density itself.

\begin{figure}[!htb]
\centering\includegraphics[width=0.86\textwidth]{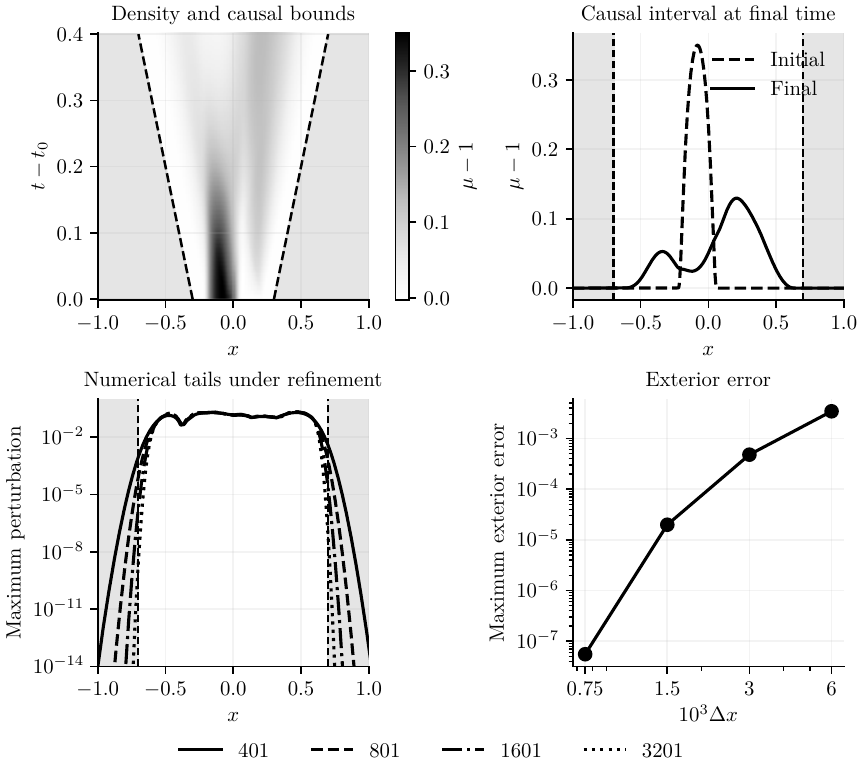}
\caption{Compact causal support of a perturbation for a nonaffine pressure. Dashed lines bound the causal domain; gray regions lie outside it. The final-time sections display the perturbation and its numerical tails on successively finer grids. The exterior error decreases under refinement; the discrete tails are not a support property of the continuous model.}
\label{fig:causal-nonaffine}
\end{figure}

\paragraph{Propagation in the numerical experiment.}
We use centered differences with Lax--Friedrichs viscosity of speed one, third-order strong-stability-preserving Runge--Kutta time stepping, and $\Delta t\le0.3\Delta x$. Figure~\ref{fig:causal-nonaffine} shows the solution through $t-t_0=0.4$ and the bounding rays $x=\pm(0.3+t-t_0)$. Artificial viscosity gives a small numerical tail outside these rays. On $401$, $801$, $1601$, and $3201$ grid points in $[-1.2,1.2]$, its maximum norm over the four first-order variables at the final time is respectively $3.43\times10^{-3}$, $4.79\times10^{-4}$, $1.99\times10^{-5}$, and $5.55\times10^{-8}$. The error in the discrete constraint is below $3\times10^{-14}$. The $L^1$ differences between solutions on successive grids decrease from $1.81\times10^{-2}$ to $6.69\times10^{-3}$. These finite tests illustrate the characteristic speed bound; they are not a proof of convergence of the numerical scheme.

\subsection{Comparison under a Lorentz boost}

\paragraph{Initial data on the boosted time slice.}
We compare evolutions in two frames related by a Lorentz boost. After transforming to the same frame, the discrepancy decreases under mesh refinement for the Lorentz-invariant viscosity and persists for the spatial Laplacian.
We compare two independent evolutions in frames related by a boost of speed $\beta=0.3$:
\be
t'=\gamma_\beta(t-\beta x),\qquad x'=\gamma_\beta(x-\beta t),\qquad
\gamma_\beta\coloneqq(1-\beta^2)^{-1/2}.
\ee
The first evolution uses the Lorentz-invariant viscosity in~\eqref{eq:conservative-tau-model-main}; the second uses artificial spatial viscosity, with a time-dependent coefficient multiplying the Laplacian of the same velocity,
\be
\del_\alpha T^{\alpha\nu}=\eps t\del_x^2U^\nu.
\ee
In the primed frame each model is solved with its own formula, using $\tau'=\sqrt{t'^2-x'^2}$ for the first and $t'\del_{x'}^2$ for the second. Both act on the unit velocity, with $\eps=0.03$ and the nonaffine pressure above.

For each model, a reference solution is computed on $2\le t\le3.6$, $|x|\le1.7$, from
\be
\mu(2,x)=1+0.12\psi((x+0.08)/0.15),\qquad
\theta(2,x)=0.04\psi((x-0.08)/0.15).
\ee
For Lorentz-invariant viscosity we also prescribe $w_t=0$, with $w=\sinh\theta$ and $w_x$ initialized by differentiation. We extract the reference solution on the tilted line corresponding to $t'=2.8$ and transform its trace to initialize a second evolution on $-1.8\le x'\le0.2$. The density is a scalar and $\theta'=\theta-\operatorname{arctanh}\beta$. For Lorentz-invariant viscosity its derivative data are transformed as well:
\be
w'=\gamma_\beta(w-\beta\sqrt{1+w^2}),\quad
w'_{t'}=\gamma_\beta^2(1-\beta v)(w_t+\beta w_x),\quad
w'_{x'}=\gamma_\beta^2(1-\beta v)(w_x+\beta w_t).
\ee
Here $v=w/\sqrt{1+w^2}$, and prime marks the observer rather than a derivative. The Laplacian model requires only the extracted density and rapidity as initial data. Thus the two frames use traces of the same reference evolution, rather than unrelated data on two horizontal time slices.

\paragraph{Covariance under mesh refinement.}
At $t'=3.1$, we compare the second evolution with the transformed reference at the same events, using bilinear interpolation of the stored reference solution. The measured error is
\be
E_h\coloneqq\int_{-1.4}^{-0.3}
\left(|\mu_{\rm evolved}-\mu_{\rm transformed}|
+|\theta_{\rm evolved}-\theta_{\rm transformed}|\right)\,dx'.
\ee
Both evolutions use centered differences, a vanishing Lax--Friedrichs term and second-order Runge--Kutta time stepping, with the parabolic time-step restriction imposed for the Laplacian. Boundary values are held at their initial traces, and errors are measured in the interior interval above.
\begin{itemize}
\item \textbf{Mesh refinement.} For $201,401,801,1601$ points in each frame, the errors for Lorentz-invariant viscosity are $3.22\times10^{-4}$, $1.84\times10^{-4}$, $9.96\times10^{-5}$, and $5.20\times10^{-5}$. The Laplacian errors are $9.16\times10^{-4}$, $9.15\times10^{-4}$, $9.22\times10^{-4}$, and $9.28\times10^{-4}$.
\item \textbf{Sensitivity checks.} An enlarged-domain check with $801$ points gives $1.01\times10^{-4}$ and $9.16\times10^{-4}$, respectively. At $801$ points, halving the evolution time step and reducing the maximum reference-output spacing from approximately $8\times10^{-4}$ to $2\times10^{-4}$ changes $E_h$ by less than $5\times10^{-9}$ for both models.
\item \textbf{Zero-boost restart.} A zero-boost restart check with $401$ points, measured on $[-0.6,0.6]$, gives $3.54\times10^{-5}$ for Lorentz-invariant viscosity and $3.41\times10^{-5}$ for spatial diffusion. This last check measures restart and interpolation errors rather than a discrepancy under the change of frame.
\end{itemize}

\begin{figure}[!htb]
\centering
\includegraphics[width=0.74\textwidth]{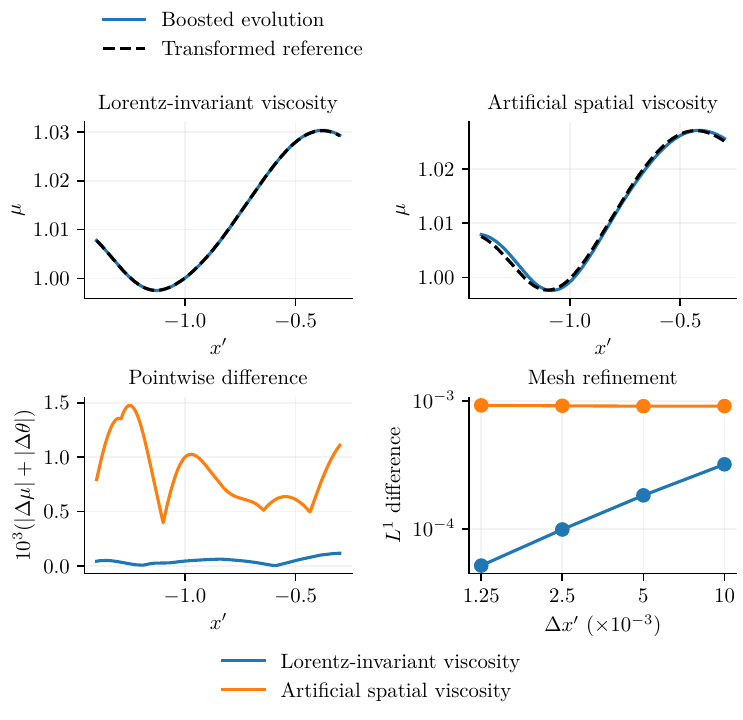}
\caption{Independent evolutions in a boosted frame compared with the transformed reference solution at $t'=3.1$. The upper panels show density on the finest mesh for Lorentz-invariant viscosity and for artificial viscosity given by a spatial Laplacian in the frame of computation; the lower panels show the combined density--rapidity difference and its $L^1$ norm under refinement. Initial traces and, for Lorentz-invariant viscosity, velocity derivative data are transformed from the reference evolution.}
\label{fig:boosted-solutions}
\end{figure}

Figure~\ref{fig:boosted-solutions} compares the profiles at the same events and displays the discrepancies on the tested meshes.

\subsection*{Acknowledgments}

 The first author was supported by the French grant ANR-23-CE40-0010 and by the European grant HORIZON-MSCA-2022-SE/101131233.  \emph{Note on AI use:} After we had written a first version of this paper, we used ChatGPT 6 to improve its presentation and produce the figures and numerical results. This was followed by a human-only review to ensure mathematical validity and clarity.

\FloatBarrier

\end{document}